\documentclass[11pt]{amsart}
\usepackage{amsmath,amsthm,amsfonts,amssymb,mathrsfs}
\usepackage{enumerate}
\usepackage[colorlinks=true, linkcolor=blue, citecolor=magenta, menucolor=black]{hyperref}
\usepackage{geometry}
\usepackage[utf8]{inputenc}
\usepackage[demo]{graphicx}
\usepackage[up]{caption}
\usepackage{subcaption}
\usepackage{pgf,tikz}
\usepackage[toc,page]{appendix}
\usetikzlibrary{arrows}
\usetikzlibrary{patterns}
\usepackage{color}

\usepackage{tikz-cd}

\makeatletter
\@namedef{subjclassname@2020}{%
  \textup{2020} Mathematics Subject Classification}
\makeatother

\numberwithin{equation}{section}
\theoremstyle{plain}
\newtheorem{thm}{Theorem}[section]

\newtheorem{lem}[thm]{Lemma}
\newtheorem{prop}[thm]{Proposition}

\newtheorem{letterthm}{Theorem}

\theoremstyle{definition}
\newtheorem{defn}[thm]{Definition}
\newtheorem*{defn*}{Definition}

\newtheorem{rem}[thm]{Remark}

\newtheorem*{example}{Example}

\newtheorem*{CRC}{Connes'\ rigidity conjecture}

\newcommand{\C}{\mathbb{C}}

\newcommand{\N}{\mathbb{N}}

\newcommand{\R}{\mathbb{R}}

\newcommand{\Z}{\mathbb{Z}}

\newcommand{\bfG}{\mathbf{G}}
\newcommand{\bfH}{\mathbf{H}}
\newcommand{\bfP}{\mathbf{P}}
\newcommand{\bfQ}{\mathbf{Q}}

\newcommand{\bfT}{\mathbf{T}}
\newcommand{\bfV}{\mathbf{V}}

\newcommand{\cC}{\mathscr{C}}

\newcommand{\cS}{\mathscr{S}}

\newcommand{\cU}{\mathscr{U}}

\newcommand{\cZ}{\mathscr{Z}}

\newcommand{\talpha}{{\widetilde{\alpha}}}

\newcommand{\tA}{{\widetilde{A}}}

\newcommand{\tE}{{\widetilde{E}}}

\newcommand{\tM}{{\widetilde{M}}}
\newcommand{\tN}{{\widetilde{N}}}

\newcommand{\oV}{{\overline V}}
\newcommand{\oU}{{\overline U}}
\newcommand{\oP}{{\overline P}}
\newcommand{\oQ}{{\overline Q}}
\newcommand{\oR}{{\overline R}}
\newcommand{\ov}{{\overline v}}
\newcommand{\ou}{{\overline u}}
\newcommand{\obfV}{\overline{\mathbf{V}}}
\newcommand{\obfP}{\overline{\mathbf{P}}}

\newcommand{\rB}{\operatorname{B}}

\newcommand{\tr}{\operatorname{tr}}

\newcommand{\Stab}{\operatorname{Stab}}
\newcommand{\Fix}{\operatorname{Fix}}

\newcommand{\inn}{\operatorname{inn}}
\newcommand{\ot}{\otimes}
\newcommand{\ovt}{\mathbin{\overline{\otimes}}}
\newcommand{\Aut}{\operatorname{Aut}}

\newcommand{\Ind}{\operatorname{Ind}} 
 
\newcommand{\id}{\operatorname{id}}
\newcommand{\GL}{\operatorname{GL}}

\newcommand{\EL}{\operatorname{EL}}

\newcommand{\SL}{\operatorname{SL}}

\newcommand{\Map}{\operatorname{Map}}

\newcommand{\Prob}{\operatorname{Prob}}
\newcommand{\supp}{\operatorname{supp}}
\newcommand{\Sub}{\operatorname{Sub}}
\newcommand{\rk}{\operatorname{rk}}
\newcommand{\bary}{\operatorname{Bar}}

\newcommand{\Rad}{\operatorname{Rad}}

\newcommand{\Gr}{\operatorname{Gr}}

\newcommand{\Sp}{\operatorname{Sp}}
\newcommand{\Cl}{\operatorname{Cl}}

\newcommand{\Out}{\operatorname{Out}}

\newcommand{\dpr}{^{\prime\prime}}

\newcommand{\eps}{\varepsilon}
\newcommand{\actson}{\curvearrowright}

\newcommand{\rC}{\operatorname{C}}
\newcommand{\rL}{\operatorname{ L}}

\newcommand{\Char}{\operatorname{Char}}

\allowdisplaybreaks

\begin{document}

\title[Charmenability of higher rank arithmetic groups]{Charmenability of higher rank arithmetic groups}

\author{Uri Bader}
\address{Faculty of Mathematics and Computer Science \\ The Weizmann Institute of Science \\ 234 Herzl Street \\ Rehovot 7610001 \\ ISRAEL}
\email{bader@weizmann.ac.il}
\thanks{UB is supported by ISF Moked 713510 grant number 2919/19}

\author{R\'emi Boutonnet}
\address{Institut de Math\'ematiques de Bordeaux \\ CNRS \\ Universit\'e Bordeaux I \\ 33405 Talence \\ FRANCE}
\email{remi.boutonnet@math.u-bordeaux.fr}
\thanks{RB is supported by ANR grant AODynG 19-CE40-0008}

\author{Cyril Houdayer}
\address{Universit\'e Paris-Saclay \\ Institut Universitaire de France \\ CNRS \\ Laboratoire de math\'ematiques d'Orsay \\ 91405 \\ Orsay  \\ FRANCE}
\email{cyril.houdayer@universite-paris-saclay.fr}
\thanks{CH is supported by Institut Universitaire de France}


\begin{abstract}
We complete the study of characters on higher rank semisimple lattices initiated in \cite{BH19,BBHP20}, the missing case being the case of lattices in higher rank simple algebraic groups in arbitrary characteristics. More precisely, we investigate dynamical properties of the conjugation action of such lattices on their space of positive definite functions. Our main results deal with the existence and the classification of characters from which we derive applications to topological dynamics, ergodic theory, unitary representations and operator algebras. Our key theorem is an extension of the noncommutative Nevo-Zimmer structure theorem obtained in \cite{BH19} to the case of simple algebraic groups defined over arbitrary local fields. We also deduce a noncommutative analogue of Margulis'\! factor theorem for von Neumann subalgebras of the noncommutative Poisson boundary of higher rank arithmetic groups.
\end{abstract}

\maketitle

%%%%%%%%%%
\section{Introduction and statements of the main results}

For any countable discrete group $\Gamma$, we denote by $\mathscr P(\Gamma) \subset \ell^\infty(\Gamma)$ the weak-$*$ compact convex set of all positive definite functions $\varphi : \Gamma \to \C$ normalized so that $\varphi(e) = 1$. We endow the space $\mathscr P(\Gamma)$ with the affine conjugation action $\Gamma \curvearrowright \mathscr P(\Gamma)$. Recall that for any positive definite function $\varphi \in \mathscr P(\Gamma)$, there is a unique (up to unitary conjugation) GNS triple $(\pi_\varphi, H_\varphi, \xi_\varphi)$, where $\pi_\varphi : \Gamma \to \mathscr U(H_\varphi)$ is a unitary representation and $\xi_\varphi \in H_\varphi$ is a unit cyclic vector such that $\varphi(\gamma) = \langle \pi_\varphi(\gamma)\xi_\varphi, \xi_\varphi\rangle$ for every $\gamma \in \Gamma$. Thus, one can always regard any positive definite function as a coefficient of a unitary representation. For standard facts on operator algebras ($\rC^*$-algebras and von Neumann algebras), we refer the reader to \cite{Ta02}.

A conjugation invariant positive definite function $\varphi \in \mathscr P(\Gamma)$ is called a {\em character} and we denote by $\Char(\Gamma) \subset \mathscr P(\Gamma)$ the weak-$\ast$ closed convex subset of all characters on $\Gamma$. Any non-trivial group $\Gamma$ always possesses at least two characters: the trivial character $1_\Gamma$ and the regular character $\delta_e$ whose GNS representation coincides with the left regular representation $\lambda : \Gamma \to \mathscr U(\ell^2(\Gamma))$. Any finite dimensional unitary representation $\pi : \Gamma \to \mathscr U(n)$ gives rise to the character $\tr_n \circ \pi$. Also, any probability measure preserving (pmp) action $\Gamma \curvearrowright (X, \nu)$ gives rise to the character $\varphi_\nu : \Gamma \to \C : \gamma \mapsto \nu(\Fix(\gamma))$ where $\Fix(\gamma) = \{x \in X \mid \gamma x = x\}$ for every $\gamma \in \Gamma$. Note that $\varphi_\nu = \delta_e$ if and only if the pmp action $\Gamma \curvearrowright (X, \nu)$ is essentially free, that is, for almost every $x \in X$, we have $\Stab_\Gamma(x) = \{e\}$.

We now review the notions of charmenable and charfinite groups that were recently introduced in \cite{BBHP20}. We denote by $\Rad(\Gamma)$ the amenable radical of $\Gamma$, that is, the largest normal amenable subgroup of $\Gamma$.

\begin{defn*}[{\cite{BBHP20}}]\label{def:charmenable}  We say that $\Gamma$ is {\em charmenable} if it satisfies the following two properties:
\begin{itemize}
\item [(P1)] Every nonempty $\Gamma$-invariant weak-$\ast$ compact convex subset $\mathscr C \subset \mathscr P(\Gamma)$ contains a fixed point, that is, a character.
\item [(P2)] Every extremal character $\varphi \in \Char(\Gamma)$ is either supported on $\Rad(\Gamma)$ or its GNS von Neumann algebra $\pi_\varphi(\Gamma)\dpr$ is amenable.
\end{itemize}
Moreover, we say that  $\Gamma$ is {\em charfinite} if it also satisfies the following three properties:
\begin{itemize} 
\item [(P3)] $\Rad(\Gamma)$ is finite.
\item [(P4)] $\Gamma$ has a finite number of isomorphism classes of unitary representations in each given finite dimension.
\item [(P5)] Every extremal character $\varphi \in \Char(\Gamma)$ is either supported on $\Rad(\Gamma)$ or its GNS von Neumann algebra $\pi_\varphi(\Gamma)\dpr$ is finite dimensional.
\end{itemize}
\end{defn*}
Recall that a tracial von Neumann algebra $M \subset \rB(\rL^2(M))$ is {\em amenable} if there exists a state $\varphi \in \rB(\rL^2(M))^*$ such that $\varphi(xT) = \varphi(Tx)$ for every $T \in \rB(\rL^2(M))$ and every $x \in M$ and $\varphi|_M$ is a faithful normal tracial state on $M$.

Before reviewing previous works and providing examples, let us first recall some striking properties of charmenable and charfinite groups. In that respect, denote by $\Sub(\Gamma) \subset 2^\Gamma$ the closed subset of all subgroups of $\Gamma$ endowed with the conjugation action $\Gamma \curvearrowright \Sub(\Gamma)$. The topology on $\Sub(\Gamma)$ induced from the product topology on $2^\Gamma$ is called the Chabauty topology. Following \cite{AGV12}, an {\em Invariant Random Subgroup}, or IRS for short, is a $\Gamma$-invariant Borel probability measure on $\Sub(\Gamma)$. Following \cite{GW14}, a {\em Uniformly Recurrent Subgroup}, or URS for short, is a nonempty $\Gamma$-invariant minimal closed subset of $\Sub(\Gamma)$.

Any charmenable group $\Gamma$ enjoys the following properties:  any normal subgroup $N \lhd \Gamma$ is either amenable or coamenable;  for any URS $X \subset \Sub(\Gamma)$, either all the elements of $X$ are contained in $\Rad(\Gamma)$ or $X$ carries an IRS;  any nonamenable unitary representation $\pi : \Gamma \to \mathscr U(H_\pi)$ weakly contains the left regular representation $\lambda$. Moreover, in that case, the $\rC^*$-algebra $\rC^*_\pi(\Gamma) \subset \rB(H_\pi)$ generated by $\pi(\Gamma)$ has a unique trace and a unique maximal proper ideal.

Any charfinite group enjoys the following properties: any normal subgroup $N \lhd \Gamma$ is either finite or has finite index;  any ergodic IRS and any URS of $\Gamma$ is finite. Any charmenable group $\Gamma$ with property (T) and such that $\Rad(\Gamma)$ is finite is charfinite. For all these facts, we refer the reader to \cite[Section 3]{BBHP20}. In particular, charfinite groups satisfy the conclusion of Margulis'\! normal subgroup theorem \cite{Ma91} and Stuck-Zimmer's stabilizer rigidity theorem \cite{SZ92}. The motivation for studying charmenable and charfinite groups comes from these fundamental results as well as other results regarding the classification of characters (see \cite{Be06, PT13, CP13, Pe14}).

The first class of charfinite groups were obtained in \cite{BH19}. More precisely, the main results of \cite{BH19} show that lattices in higher rank connected simple Lie groups with finite center are charfinite. New classes of charmenable and charfinite groups were subsequently obtained in \cite{BBHP20}. Indeed, \cite[Theorem A]{BBHP20} shows that irreducible lattices in products of (at least two) simple algebraic groups are charmenable (resp.\! charfinite if one of the factors has property (T)). Morever, for every $d \geq 2$ and every nonempty (possibly infinite) set of primes $S \subset \mathscr P$, the $S$-arithmetic group $\SL_d(\Z[S^{-1}])$ is charfinite. We refer the reader to the recent survey \cite{Ho21} for further details and results.

The goal of the present paper is to complete the above picture by showing that lattices in higher rank simple algebraic groups defined over a local field $k$ of arbitrary characteristic are charfinite.

\begin{letterthm} \label{thm:AG-simple}
Let $k$ be a local field. Let $\mathbf G$ be an almost $k$-simple connected algebraic $k$-group such that $\rk_k(\mathbf G) \geq 2$. Then every lattice $\Gamma < \bfG(k)$ is charfinite.
\end{letterthm}

Theorem \ref{thm:AG-simple} provides several new classes of charfinite groups. 

\begin{example}
Let $d \geq 3$ be an integer and $k$ a local field. Then any lattice $\Gamma < \SL_d(k)$ is charfinite. In particular, let $p \in \mathscr P$ be a prime and $q = p^r$ for $r \geq 1$. Denote by $\bfQ_p$ the local field of $p$-adic numbers and by $\mathbf F_q((t))$ the local field of formal power series in one variable $t$ over the finite field $\mathbf F_q$. Then any lattice $\Gamma < \SL_d(\bfQ_p)$ and any lattice $\Gamma < \SL_d(\mathbf F_q((t)))$ is charfinite. In particular, the lattice $\SL_d(\mathbf F_q[t^{-1}]) < \SL_d(\mathbf F_q((t)))$ is charfinite. We also refer to \cite{LL20} for the character classification of $\SL_d(\mathbf F_q[t^{-1}])$.
\end{example}

Before stating our next result regarding charmenability of higher rank arithmetic groups, let us review some terminology. Let $K$ be a global field and ${\bfG}$ an almost $K$-simple connected algebraic $K$-group. Let $S$ be a (possibly empty, possibly infinite) set of non-Archimedean inequivalent absolute values on $K$. Let $\mathscr{O}<K$ be the ring of integers and $\mathscr{O}_S$ the corresponding localization, that is,
\[ \mathscr{O}_S = \{ \alpha\in K \mid \forall s\in S, s(\alpha) \leq 1\}.  \]
Fix an injective $K$-representation $\rho:{\bfG} \to \GL_n$ and write
\[ \Lambda_S=\rho^{-1}(\GL_n(\mathscr{O}_S))\le {\bf G}(K). \]
The triple $(K,{\bf G},S)$ is said to be of
\begin{itemize}
\item \emph{compact type}
if for every absolute value $v$ on $K$, the image of $\Lambda_S$ in ${\bf G}(K_v)$ is bounded;
\item \emph{simple type}
if there exists a unique absolute value $v$ on $K$ such that the image of $\Lambda_S$ in ${\bf G}(K_v)$ is unbounded;
\item \emph{product type} otherwise.
\end{itemize}
The triple $(K,{\bf G},S)$ is said to be \emph{of higher rank} if it is either of product type or of simple type and $\rk_{K_v}({\bf G}) \geq 2$. A subgroup $\Gamma\le {\bf G}(K)$ is called \emph{$S$-arithmetic} if it is commensurable with $\Lambda_S$. It is called \emph{arithmetic} if it is $S$-arithmetic for some $S$ as above and we regard its type as the type of $(K,{\bf G},S)$.

Combining Theorem \ref{thm:AG-simple} with \cite[Theorem A]{BBHP20}, we infer that all higher rank arithmetic groups are charmenable.

\begin{letterthm} \label{thm:AG-general}
Let $K$ be a global field and $\mathbf G$ an almost $K$-simple connected algebraic $K$-group. Let $S$ be a (possibly empty, possibly infinite) set of non-Archimedean inequivalent absolute values on $K$. Then any higher rank $S$-arithmetic group $\Gamma \leq \mathbf G(K)$ is charmenable.

Assume further that there exists an absolute value $\upsilon$ on $K$ such that $\mathbf G(K_\upsilon)$ has property {\em (T)} and for which the image of $\Gamma$ in $\mathbf G(K_\upsilon)$ is unbounded. If either $S$ is finite or $\mathbf G$ is simply connected, then $\Gamma$ is charfinite.
\end{letterthm}

For a $\rC^*$-algebra $A\subset \mathrm B(H)$ and $n \geq 1$, $\mathrm M_n(A) = \mathrm M_n(\C) \otimes A \subset \mathrm B( H^{\oplus n})$ is naturally a $\rC^*$-algebra. Let $A, B$ be $\rC^*$-algebras. A linear map $\Phi : A \to B$ is said to be \emph{unital completely positive} (ucp) if $\Phi$ is unital and if for every $n \geq 1$, the linear map $\Phi^{(n)} : \mathrm M_n(A) \to \mathrm M_n(B) : [a_{ij}]_{ij} \mapsto [\Phi(a_{ij})]_{ij}$ is positive. Any unital $\ast$-homomorphism $\pi : A \to B$ is a ucp map. When $A$ or $B$ is commutative, any unital positive linear map $\Phi : A \to B$ is automatically ucp (see e.g.\! \cite[Theorems 3.9 and 3.11]{Pa02}). 

As in \cite{BH19}, our proof of Theorem \ref{thm:AG-simple} is based on a noncommutative analogue of Nevo-Zimmer structure theorem (see \cite{NZ97,NZ00}) which has further applications. 

\begin{letterthm}\label{thm:NCNZ}
Let $k$ be a local field. Let $\mathbf G$ be an almost $k$-simple connected algebraic $k$-group such that $\rk_k(\mathbf G) \geq 2$ and set $G = \mathbf G(k)$. Let $\mathbf P < \mathbf G$ be a minimal parabolic $k$-subgroup and set $P = \mathbf P(k)$. Let $M$ be an ergodic $G$-von Neumann algebra and $E : M \to L^\infty(G/P)$ a faithful normal ucp $G$-map. The following dichotomy holds:
\begin{itemize}
\item Either $E$ is $G$-invariant, that is, $E(M) = \C 1$.
\item Or there exist a proper parabolic $k$-subgroup $\mathbf P < \mathbf Q < \mathbf G$ and a $G$-equivariant unital normal embedding $\iota : L^\infty(G/Q) \hookrightarrow M$ where $Q = \mathbf Q(k)$ such that $E \circ \iota : L^\infty(G/Q) \hookrightarrow L^\infty(G/P)$ is the canonical unital normal embedding.
\end{itemize}
\end{letterthm}

Both the statement and the proof of Theorem \ref{thm:NCNZ} are similar to \cite[Theorem 5.1]{BH19}, but extra difficulties appear in the proof to handle simple algebraic groups in positive characteristic. Assume that $E : M \to L^\infty(G/P)$ is not invariant. Firstly, we construct in Theorem \ref{First half} an abelian $G$-von Neumann subalgebra $M_0 \subset M$ for which the restriction $E|_{M_0} : M_0 \to L^\infty(G/P)$ is not invariant either. We point out that compared to the proof of \cite[Theorem 5.1]{BH19}, some of the arguments have been simplified by making use of the so-called {\em maximal compact models}. We are grateful to Amine Marrakchi for sharing this idea with us. Secondly, we prove in Theorem \ref{thm:gauss-stat} a generalization of the commutative Nevo-Zimmer structure theorem for stationary actions of higher rank simple algebraic groups on standard probability spaces (see \cite[Theorem 1]{NZ00} for actions of higher rank connected simple Lie groups). Along the way, we also prove some useful facts regarding closed subgroups of algebraic groups that are of independent interest. Then Theorem \ref{thm:AG-simple} follows from Theorem \ref{thm:NCNZ} by using an induction argument in a similar fashion as \cite{BBHP20}.

For a countable discrete group $\Lambda$ and a nonsingular action $\Lambda \curvearrowright (Y, \eta)$ on a standard probability space, we denote by $L(\Lambda \curvearrowright Y)$ the {\em group measure space von Neumann algebra}. The von Neumann algebra $L(\Lambda \curvearrowright Y)$ is generated by a copy of the {\em group von Neumann algebra} $L(\Lambda) = \{u_\gamma \mid \gamma \in \Lambda\}\dpr \subset \rB(\ell^2(\Lambda))$ and a copy of $L^\infty(Y, \eta)$ in such a way that
$$\forall \gamma \in \Lambda, \forall F \in L^\infty(Y, \eta), \quad u_\gamma F u_\gamma^* = F \circ \gamma^{-1}.$$
If $\Lambda \curvearrowright (Y, \eta)$ is essentially free and ergodic, then $L(\Lambda \curvearrowright Y)$ is a von Neumann factor whose type coincides with the type of the action $\Lambda \curvearrowright (Y, \eta)$ (see e.g.\! \cite[Theorem XIII.1.7]{Ta03}). For lattices in higher rank simple algebraic groups, we present yet another application of Theorem \ref{thm:NCNZ} (or rather Theorem \ref{thm:NCNZ-lattices}). Our next result may be regarded as a noncommutative analogue of Margulis'\! celebrated factor theorem (see \cite[Theorem IV.2.11]{Ma91}).

\begin{letterthm}\label{thm:NCFT}
Let $k$ be a local field. Let $\mathbf G$ be an almost $k$-simple connected algebraic $k$-group such that $\rk_k(\mathbf G) \geq 2$ and set $G = \mathbf G(k)$. Let $\mathbf P < \mathbf G$ be a minimal parabolic $k$-subgroup and set $P = \mathbf P(k)$. Assume that $\mathbf G$ is center free, that is, $\mathscr Z(\bfG) = \{e\}$. 

Then for every intermediate von Neumann subalgebra $L(\Gamma) \subset M \subset L(\Gamma \curvearrowright G/P)$, there exists a unique parabolic $k$-subgroup $\mathbf P < \mathbf Q < \mathbf G$ such that $ M = L(\Gamma \curvearrowright G/Q)$, where $Q = \mathbf Q(k)$.

In particular, the $k$-rank $\rk_k(\bfG)$ is an invariant of the inclusion $L(\Gamma) \subset L(\Gamma \curvearrowright G/P)$.
\end{letterthm}

The group measure space von Neumann algebra $L(\Gamma \curvearrowright G/P)$ is an amenable type ${\rm III}_1$ factor. In case $k = \R$, $L(\Gamma \curvearrowright G/P)$ coincides with the noncommutative Poisson boundary of the lattice $\Gamma$ (see \cite{Fu67, Iz04}).  Theorem \ref{thm:NCFT} was recently announced and stated in \cite{Ho21} in the case of lattices in higher rank simple Lie groups with trivial center (see \cite[Corollary F]{Ho21}).

Besides this introduction and the next preliminary section, the paper contains three other sections. In Section \ref{reduce to commutative}, we prove the first half of Theorem \ref{thm:NCNZ} by constructing a non-trivial $G$-invariant abelian von Neumann subalgebra $M_0 = L^\infty(X, \nu)$ for which the corresponding nonsingular action $G \curvearrowright (X, \nu)$ has large stabilizers. In Section \ref{Gauss}, we prove the second half of Theorem \ref{thm:NCNZ}  and we generalize the Gauss map technique introduced by Nevo-Zimmer in \cite{NZ00} to algebraic groups defined over a local field of arbitrary characteristic. In Section \ref{main results}, we prove Theorems \ref{thm:AG-simple}, \ref{thm:AG-general} and \ref{thm:NCFT}. We also discuss the relevance of Theorem \ref{thm:NCFT} regarding Connes'\! rigidity conjecture for group von Neumann algebras of higher rank lattices.

\subsection*{Acknowledgments}

We thank the anonymous referee for carefully reading our paper and for providing useful comments and remarks that helped improve the exposition of the paper.

{
  \hypersetup{linkcolor=black}
  \tableofcontents
}

%%%%%%%%%%
\section{Preliminaries}

%%%%%
\subsection{Structure of (almost) simple algebraic $k$-groups}
\label{AGnot}

Let $k$ be a local field, that is, $k$ is a nondiscrete locally compact field of arbitrary characteristic. An {\em algebraic} $k$-{\em group} $\bfG$ is an algebraic group that is defined over $k$. In this paper, all algebraic $k$-groups $\bfG$ are assumed to be affine (or linear). We denote by $\bfG^0$ the (Zariski) connected component in $\bfG$ of the identity element. Then $\bfG$ is (Zariski) connected if $\bfG = \bfG^0$.

Let $\bfG$ be a connected algebraic $k$-group. We say that $\bfG$ is {\em semisimple} (resp.\! {\em reductive}) if its radical (resp.\! unipotent radical) is trivial. We say that $\bfG$ is {\em absolutely almost simple} if it is semisimple and all its proper algebraic normal subgroups are finite. We say that $\bfG$ is {\em almost $k$-simple} if it is semisimple and all its proper $k$-closed normal subgroups are finite. Note that since $\bfG$ is connected, finite normal subgroups of $\bfG$ are contained in $\mathscr Z(\bfG)$.

Let $\bfG$ be an almost $k$-simple connected algebraic $k$-group. We will be using the following notation throughout:
\begin{itemize}
\item $\bfT$ will denote a maximal $k$-split torus of $\bfG$, $\Phi^+$ a choice of positive roots, and $\Delta$ the corresponding set of simple positive roots.
\item $\bfP$ will denote the corresponding minimal parabolic $k$-subgroup, and more generally, for every subset $\theta \subset \Delta$, $\bfT_\theta \subset \bfT$ and $\bfP_\theta$ will denote the corresponding torus and parabolic $k$-subgroup.
\item $\bfV$ (resp.\! $\bfV_\theta$) will denote the unipotent radical of $\bfP$ (resp.\! $\bfP_\theta$). For each $\theta$, one has the Levi decomposition $\bfP_\theta = \bfH_\theta \bfV_\theta$, where the reductive group $\bfH_\theta$ is the centralizer of $\bfT_\theta$ inside $\bfG$.
\item The opposite parabolic $k$-subgroups will be denoted by $\obfP$, $\overline\bfP_\theta$ and their unipotent radical will be naturally denoted by $\obfV$, $\obfV_\theta$.
\item We will consider the corresponding groups of $k$-points: $G = \bfG(k)$, $T= \bfT(k)$, $P = \bfP(k)$, $V = \bfV(k)$, $\oP = \obfP(k)$, $\oV = \obfV(k)$, and similary, $T_\theta$, $P_\theta$, $V_\theta$, $H_\theta$, $\oP_\theta$, $\oV_\theta$.
\end{itemize}

With this notation, for every $\theta$, the product map $\obfV_\theta \times \bfP_\theta \to \bfG$ is a $k$-isomorphism onto a Zariski dense and open subset of $\bfG$. At the level of $k$-points, the set $\oV_\theta P_\theta$ is co-null inside $G$ with respect to the Haar measure. This implies the following classical observation.

\begin{lem}\label{LUdec}
There is a measure class preserving isomorphism $G \simeq \oV_\theta \times P_\theta$, where each of the groups is endowed with its Haar measure class. This isomorphism is equivariant under the left-right action of $\oV_\theta \times P_\theta$. Since $H_\theta$ normalizes $\oV_\theta$ and is contained in $P_\theta$ this isomorphism maps the left action $H_\theta \actson G$ to the product action 
\[H_\theta \times \oV_\theta \times P_\theta \to \oV_\theta \times P_\theta : (h,v,p) \mapsto (hvh^{-1},hp).\]
\end{lem}

We will also denote by $G^+$ the subgroup of $G$ generated by $R_u({\bf P})(k)$ 
where ${\bf P}$ runs through the set of all proper parabolic $k$-subgroups of ${\bf G}$ and $R_u({\bf P})$ denotes the unipotent radical of 
${\bf P}$, see \cite[I.1.5.2]{Ma91}.
Note that by \cite[Theorem I.2.3.1]{Ma91}, $G^+$ is a closed, normal and cocompact subgroup of $G$.
When ${\bf G}$ is $k$-anisotropic, $G^+$ is trivial. When ${\bf G}$ is $k$-isotropic, if $k$ is of characteristic $0$ then $G/G^+$is finite,
if $k=\mathbb{R}$ then $G^+$ is the identity component in $G$ and if $k=\mathbb{C}$ then $G^+=G$. 

%%%%%
\subsection{Group actions on operator algebras}

We will be interested in group actions on $\rC^*$-algebras and von Neumann algebras. We refer to \cite[Section 2.1]{BH19} for the precise continuity requirements about such actions, examples, and connections between group actions on $\rC^*$-algebras and von Neumann algebras. 

\begin{defn}
Consider a locally compact second countable (lcsc) group $G$ with a closed subgroup $P < G$ and a von Neumann action $\sigma : P \actson M$. We have two commuting actions of $G$ and $P$ on the tensor product von Neumann algebra $L^\infty(G) \ovt M$; the $G$-action is given by the automorphisms $\lambda_g \ot \id$ for all $g \in G$, while the $P$-action is a product action described by the automorphisms $\rho_p \ot \sigma_p$, $p \in P$. Here $\lambda$ and $\rho$ denote the actions of $G$ on $L^\infty(G)$ induced by left and right multiplication respectively. Then $\Ind_P^G (M) =  (L^\infty(G) \ovt M)^{(\rho \otimes \sigma)(P)}$ is the {\em induced von Neumann algebra} and $\lambda \otimes \id : G \curvearrowright \Ind_P^G (M)$ is the {\em induced action}. We will simply denote by $\tM=(L^\infty(G) \ovt M)^{P}$ the induced von Neumann algebra.
\end{defn}

Given $P < G$ as above, if $N,M$ are two $P$-von Neumann algebras with a normal, $P$-equivariant unital completely positive (ucp) map $E: M \to N$, then $\tE := \id \ot E$ is a normal, $G$-equivariant ucp map $\tE: \tM \to \tN$. So, induction is functorial in this sense.

Let us introduce some more language attached to the category of $G$-von Neumann algebras with morphisms given by normal equivariant ucp maps. We will use the letter $\sigma : G \curvearrowright M$ to denote the action of $G$ on the von Neumann algebra $M$.

\begin{defn}
Consider two $G$-von Neumann algebras $N$ and $M$ and a $G$-equivariant normal ucp map $E: M \to N$. Then $E$ is said to be {\em $G$-invariant} if $E(\sigma_g(x)) = E(x)$ for every $x \in M$, $g \in G$. Equivalently, this means that $E(M) \subset N^G$. We say that $E$ is {\em faithful} if $E(x^*x) = 0$ implies $x = 0$, for every $x \in M$. The {\em support} of $E$ is the smallest projection $p \in M$ such that $E(p) = 1$. Thus $E$ is faithful if and only its support is $1$.
\end{defn}

The plan in this section is to make further observations about induction in specific situations. For that, we need to discuss compact models associated with von Neumann actions.
Given a locally compact group $G$ and a $G$-von Neumann algebra $M$, the $G$-action on $M$ need not be norm continuous in general. However, there always exists a weak-* dense C*-subalgebra $A \subset M$ which is globally invariant under $G$ and on which the action $G \curvearrowright A$ is norm continuous (see the proof of \cite[Proposition XIII.1.2]{Ta03}). By analogy with the commutative case, we call $A$ together with its $G$-action, a {\em compact model} of $M$. It is sometimes convenient to find separable compact models, but in general there is no canonical choice of such. In contrast, there is always a canonical choice of a compact model if one ignores this separability condition.

\begin{defn}
Given a $G$-von Neumann algebra $M$, the set of elements $x \in M$ such that the map $g \in G \mapsto \sigma_g(x)$ is norm continuous is a $G$-invariant C*-subalgebra of $M$. It is the largest compact model for $M$, and called the {\em maximal compact model}. Its elements are said to be {\em $G$-continuous}\footnote{There is a similar terminology in \cite[Section 5]{BBHP20}, but it refers to a different notion.}.
\end{defn}

For example, when $G$ is discrete the maximal compact model is the whole of $M$. Here is our most fundamental example: 

\begin{lem}\label{MCM0}
Consider a lcsc group $G$ with a closed subgroup $H < G$. Endow $G/H$ with its unique $G$-invariant measure class. Then the maximal compact model for the translation action $G \actson L^\infty(G/H)$ is contained in the $\rC^*$-algebra $C_b(G/H)$ of all bounded continuous functions on $G/H$.
\end{lem}

\begin{proof}
Let $f \in L^\infty(G/H)$ be $G$-continuous and take an approximate unit $f_n \in C_c(G)$, $n \in \N$. Then $f_n \ast f$ converges to $f$ in norm. Moreover, each $f_n \ast f$ is continuous on $G/H$, hence so is~$f$.
\end{proof}

\begin{rem}
As we will see, the above lemma will allow us to avoid the use of separability arguments. So the content of this section does not make any separability assumption on the involved algebras. This is not so important for our purposes, since our noncommutative Nevo-Zimmer theorem can be reduced to the version where the algebra $M$ is separable, but we point it out anyway.
\end{rem}

We will need the following refinement of the previous lemma.

\begin{lem}\label{MCM1}
Consider a lcsc group $G$ and a von Neumann algebra $N$. Then $\lambda \ot \id$ is a $G$-action on $L^\infty(G) \ovt N$. The maximal compact model of this action is contained in the $\rC^*$-algebra $C_b(G,N)$ of bounded, norm continuous functions $G \to N$ (which is naturally embedded inside $L^\infty(G) \ovt N$ as the multiplier algebra of $C_0(G) \ot_{\min} N$).
\end{lem}

\begin{proof}
Take a $G$-continuous element $f \in L^\infty(G) \ovt N$. For every $\psi \in N_*$, the map $\id \ot \psi$ is $G$-equivariant and norm continuous, so it maps $f$ to a $G$-continuous element of $L^\infty(G)$. Thanks to Lemma \ref{MCM0}, this gives a norm bounded map $N_* \to C_b(G)$. In particular, for every $g \in G$, we may compose this map with the evaluation state at $g$, to get a linear functional on $N_*$, that is, an element $f(g) \in N$. Moreover, we observe that since $f$ is $G$-continuous, the function $f: G \to N$ obtained this way is norm continuous. 

Denote by $A \subset L^\infty(G) \ovt N$ the maximal compact model for the $G$-action. By the previous paragraph, we obtain a *-homomorphism $\iota: A \to C_b(G,N)$. It is easily seen that $\iota$ is injective. Moreover $A$ contains $C_0(G) \ot_{\min} N$ and it is clear that on this algebra $\iota$ corresponds to the identification $C_0(G) \ot_{\min} N \simeq C_0(G,N)$. So, starting with $f_1 \in C_0(G) \ot_{\min} N$ and $f_2 \in A$, $\iota(f_1f_2) = \iota(f_1)\iota(f_2) \in C_0(G,N)$. Hence $f_1f_2 \in C_0(G) \ot_{\min} N$, proving that $f_2$ belongs to the desired multiplier algebra.
\end{proof}

Consider two $G$-von Neumann algebras $M$ and $N$, with a $G$-equivariant normal ucp map $E: M \to N$. Then $E$ carries any compact model of $M$ into the maximal compact model of $N$. In the special case where $N$ is of the form $L^\infty(G/H)$ for some closed subgroup $H < G$, Lemma \ref{MCM0} ensures that $E$ maps the maximal compact model $A$ of $M$ inside $C_b(G/H)$. Composing with the evaluation map at the coset $H \in G/H$, we obtain an $H$-invariant state $\phi: A \to \C$. This observation is used to prove the next proposition.

\begin{prop}\label{embedding}
Let $G$ be a lcsc group and $P < G$ a closed subgroup. If $M$ is a $G$-von Neumann algebra with a faithful normal $G$-equivariant ucp map $E: M \to L^\infty(G/P)$, then there exists a $P$-von Neumann algebra $N$ with a normal $P$-invariant state $\psi: N \to \C$ and a von Neumann embedding $\iota: M \to \tN = (L^\infty(G) \ovt N)^P$ such that $E = (\id \ot \psi) \circ \iota$.
\end{prop}

The proposition relies on the following classical lemma.

\begin{lem}
Consider two von Neumann algebras with normal states $(M,\phi)$ and $(N,\psi)$. Consider a weakly dense $*$-subalgebra $A \subset M$ and a state preserving $*$-homomorphism $\alpha: A \to N$. If $\psi$ is faithful on $\alpha(A)' \cap N$, then $\alpha$ extends to a normal $*$-homomorphism $M \to N$. 
\end{lem}
\begin{proof}
Without loss of generality, we may replace $N$ by its subalgebra $\alpha(A)\dpr$. In this case, the assumption is that $\psi$ is faithful on the center of $N$. This implies that the associated GNS representation $\pi_\psi$ is faithful on $N$. Indeed, there exists a unique central projection $z \in \mathscr Z(N)$ such that $\ker(\pi_\psi) = N z$. Then we have $\psi(z) = \langle \pi_\psi(z)\xi_\psi, \xi_\psi\rangle = 0$. Since $\psi|_{\mathscr Z(N)}$ is faithful, it follows that $z = 0$ and so $\pi_\psi$ is faithful on $N$.

Since $\alpha$ is state preserving, it induces a Hilbert space isomorphism $U: L^2(A,\phi) \to L^2(N,\psi)$. By density of $A$ inside $M$, we have the equality $L^2(A,\phi) = L^2(M,\phi)$. Then we may define a normal $*$-homomorphism $\talpha: M \to B(L^2(N,\psi))$ by the formula $\talpha(x) = U\pi_\phi(x)U^*$, for all $x \in M$. We observe that for all $x \in A$, $\talpha(x) = \pi_\psi(\alpha(x))$. In particular, $\talpha$ maps $M$ into $\pi_\psi(N)$. We may compose this morphism with $\pi_\psi^{-1}$ to get the desired extension.
\end{proof}

\begin{proof}[Proof of Proposition \ref{embedding}]
Consider $E: M \to L^\infty(G/P)$ as in the statement of the proposition. Denote by $A \subset M$ the maximal compact model for the $G$-action. By Lemma \ref{MCM0}, $E$ maps $A$ into $C_b(G/P)$, and thus composing $E|_{A}$ with the evaluation map at $P$, we obtain a $P$-invariant state $\psi$ on $A$. Denote by $N := \pi_\psi(A)\dpr$ the von Neumann algebra generated by $A$ in the GNS representation associated to $(A,\psi)$. We sill denote by $\sigma : G \curvearrowright A$ the norm continuous action.

The embedding $\iota$ is first defined on $A$ by the formula $\iota(a)(g) = \pi_\psi(\sigma_{g^{-1}}(a))$, for all $a \in A$, $g \in G$. One easily checks that this defines a $G$-equivariant $*$-homomorphism $$A \to C_b(G,N)^P \subset (L^\infty(G) \ovt N)^P.$$ 
Moreover, for every $a \in A$, $g \in G$, we have
\[(\id \ot \psi)(\iota(a))(gP) = \psi(\iota(a)(g)) = \psi(\sigma_{g^{-1}}(a)) = E(\sigma_{g^{-1}}(a))(P) = E(a)(gP).\]
So $(\id \ot \psi) \circ \iota = E$.

We now check that $\iota$ as above extends to $M$. If $\nu_P$ denotes a faithful normal state on $L^\infty(G/P)$, then $\phi_1 := \nu_P \circ E$ and $\phi_2 := \nu_P \circ (\id \ot \psi)$ are normal states on $M$ and $\tN$ respectively, and we have $\phi_2 \circ \iota = \phi_1$. So by the previous lemma, we only need to check that $\phi_2$ is faithful on $\iota(A)' \cap \tN$. Equivalently, we need to check that $\id \ot \psi$ is faithful on $\iota(A)' \cap \tN$. This follows from the next two claims:

{\bf Claim 1.} $\iota(A)' \cap \tN = (L^\infty(G) \ovt \cZ(N))^P$ (which coincides with the center of $\tN$).

Note that $\iota(A)' \cap \tN$ is globally $G$-invariant. So we only need to check that every $G$-continuous element $f \in \iota(A)' \cap \tN$ belongs to $L^\infty(G) \ovt \cZ(N)$. By Lemma \ref{MCM1}, $f$ may be viewed as a continuous function $G \to N$. So for every $a \in A$, we may view the equality $f\iota(a) = \iota(a)f$ inside $C_b(G,N)$. This gives $f(g)\pi_\psi(\sigma_{g^{-1}}(a)) = \pi_\psi(\sigma_{g^{-1}}(a))f(g)$ for every $g \in G$, $a \in A$. Since $\pi_\psi(A)$ is weakly dense in $N$, we conclude that $f(g) \in \cZ(N)$ for every $g \in G$. Hence $f \in C_b(G,\cZ(N))^P \subset (L^\infty(G) \ovt \cZ(N))^P$.

{\bf Claim 2.} $E_\psi := \id \otimes \psi$ is faithful on $(L^\infty(G) \ovt \cZ(N))^P$.

Since $\psi$ is a vector state on $N$ attached to a cyclic vector, $\psi$ is faithful on $\cZ(N)$. Hence $E_\psi = \id \ot \psi$ is faithful on $L^\infty(G) \ovt \cZ(N)$.
\end{proof}

When considering induced algebras, we will also use the following variation of Lemma \ref{MCM1}.

\begin{lem}\label{MCM2}
Consider a lcsc group $G$ with a closed cocompact subgroup $Q < G$ and consider a $Q$-von Neumann algebra $N$ with its induced $G$-von Neumann algebra $\tN = (L^\infty(G) \ovt N)^Q$. Denote by $A \subset N$ the maximal compact model for the $Q$-action. Then the maximal compact model $B \subset \tN$ for the induced $G$-action is the $\rC^*$-subalgebra $C_b(G,A)^Q$ consisting of bounded, norm continuous, $A$-valued equivariant functions.
\end{lem}
\begin{proof}
By Lemma \ref{MCM1}, $B$ is contained in $C_b(G,N)$, and thus in $C_b(G,N)^Q$. Moreover, for $f \in B$ and $g \in G$, we know that the map $h \in G \mapsto f(hg) \in N$ is norm continuous. Restricting this map to $gQg^{-1}$, we deduce that the map $q \in Q \mapsto f(gq) = \sigma_q(f(g))$ is continuous. So $f(g) \in A$, and thus $f \in C_b(G,A)^Q$, as desired. The converse inclusion is easy.
\end{proof}

%%%%%
\subsection{C*-algebraic lemmas}

In this section, we record lemmas about extremal states on C*-algebras, which will be used to make the passage from the noncommutative setting to the commutative setting explicit. Our C*-algebras will be assumed to be unital. When $A$ is a $\rC^*$-algebra, we denote by $\mathscr S(A)$ the state space of $A$.

\begin{defn}
A state on a unital C*-algebra $A$ is called {\em approximately extremal} if it belongs to the weak-* closure of the set of extremal states on $A$.
\end{defn}

\begin{lem}\label{C*1}
Consider two C*-algebras $A \subset B$. Then every extremal state on $A$ is the restriction of an approximately extremal state on $B$.
\end{lem}
\begin{proof}
Denote by $r: \cS(B) \to \cS(A)$ the restriction map. It follows from Hahn-Banach theorem that $r$ is surjective. Denote by $\cS \subset \cS(B)$ the set of approximately extremal states on $B$. We know that the closed convex hull of $r(\cS)$ is $r(\cS(B)) = \cS(A)$. Hence, by Krein-Millman theorem, the closed set $r(\cS)$ contains the extremal points in $\cS(A)$. This is what we wanted.
\end{proof}

\begin{lem}\label{C*int}
If $B$ is a C*-algebra and $A \subset B$ is a central $\rC^*$-subalgebra, then the restriction to $A$ of an extremal state of $B$ is extremal on $A$.
\end{lem}
\begin{proof}
Take an extremal state $\phi \in \cS(B)$. We only need to check that $\phi$ is multiplicative on $A$. Take $a \in A$, $0 \leq a \leq 1$. Since $a$ is in the center of $B$, the positive linear functional $\phi_a: x \in B \mapsto \phi(ax)$ is such that $0 \leq \phi_a \leq \phi$. By extremality of $\phi$, it follows that $\phi_a = \lambda \phi$ for some scalar $\lambda$. Evaluating at $1$ gives $\lambda = \phi(a)$, which proves that $a$ belongs to the multiplicative domain of $\phi$. Thus $\phi$ is multiplicative on $A$, as desired.
\end{proof}

\begin{lem}\label{C*2}
Consider a lcsc group $G$ with a closed cocompact subgroup $Q < G$ and a C*-action on a unital C*-algebra $Q \actson A$. Consider the induced $\rC^*$-algebra $B := C_b(G,A)^Q$, together with its $G$-action by left translation. For $g \in G$, denote by $e_g: B \to A$ the evaluation morphism at $g$ and consider the set $\cS := \{ \psi \circ e_g \mid \psi \in \cS(A), g \in G\} \subset \cS(B)$. Then
\begin{enumerate}[a)]
\item $\cS$ is weak-* closed inside $\cS(B)$;
\item Every extremal state on $B$ belongs to $\cS$, hence so does every approximately extremal state;
\item Every state in $\cS$ is fixed by some conjugate of the kernel of the action $Q \actson A$.
\end{enumerate}
\end{lem}
\begin{proof}
Before proving each of these statements, let us prove that $\cS$ coincides with the subset $\cS' \subset \cS(B)$ consisting of states whose restriction to the C*-subalgebra $C_b(G)^Q \simeq C(G/Q)$ is extremal. It is clear that $\cS \subset \cS'$. Conversely, take $\phi \in \cS'$. 

There exists $x \in G/Q$ such that $\phi$ coincides with the evaluation map at $x$ on $C(G/Q)$. Write $x = gQ$ for some $g \in G$.
Denote by $I \subset C_b(G,A)^Q$ the kernel of the evaluation morphism $e_g: C_b(G,A)^Q \to A$ at the point $g$. By $Q$-equivariance, $I$ coincides with the ideal of functions which vanish on the whole coset $gQ$. Note that $C(G/Q)$ is in the multiplicative domain of $\phi$. 

Fix $f \in I$ and $\eps > 0$. By continuity of $f$, we may find a function $f_0 \in C(G/Q)$ which is equal to $1$ at the point $x = gQ$, and such $\Vert f f_0\Vert < \eps$. Thus $|\phi(f)| = |\phi(f_0)\phi(f)| = |\phi(f_0f)| \leq \Vert f_0f\Vert < \eps$. Since this holds for every $\eps > 0$, we find that $\phi(f) = 0$, and thus, $\phi$ vanishes on $I$. So $\phi$ factors to a state on $A$ via the evaluation map $e_g$. This implies that $\phi \in \cS$, as desired.\\
a) Since $C(G/Q)$ is abelian, the set of its extremal states is weak-* closed, equal to the set of multiplicative characters. So $\cS = \cS'$ is weak-* closed as well.\\
b) By Lemma \ref{C*int}, the restriction of an extremal state on $B$ to the central subalgebra $C(G/Q)$ is still extremal.\\
c) Let $g \in G$, $\psi \in \cS(A)$ and $\phi = \psi \circ e_g \in \cS(B)$. Denote by $R$ the kernel of the action $Q \actson A$. Take $f \in C_b(G,A)^Q$ and $h \in gRg^{-1}$. Denoting by $r := g^{-1}h^{-1}g \in R$, we have 
\[\phi(\sigma_h(f)) = \psi(f(h^{-1}g)) = \psi(f(gr)) = \psi(\sigma_r(f(g))) = \psi(f(g)) = \phi(f).\qedhere\]
\end{proof}

%%%%%%%%%%

\section{Reduction to the commutative setting}
\label{reduce to commutative}

Our proof of Theorem \ref{thm:NCNZ} is split in two halves. 
The first half is the following theorem, which simultaneaously achieves two goals: 
\begin{itemize}
\item it reduces to the case where $M = L^\infty(X,\nu)$ is a commutative von Neumann algebra;
\item and similarly to Nevo-Zimmer's approach, it reduces to the case where the stabilizer of almost every point of $X$ has positive dimension in $G$.  
\end{itemize}
The second half will then be to deduce the conclusion of Theorem \ref{thm:NCNZ} from there. This will be achieved in Section \ref{Gauss}, by adapting the Gauss map trick from \cite{NZ00} to the general setting of algebraic groups over local fields. In fact, to be able to use this Gauss map in positive characteristics, one needs to know a bit more than just positive dimension of point stabilizers.

\begin{thm}\label{First half}
Let $k$ be a local field. Let $\mathbf G$ be an almost $k$-simple connected algebraic $k$-group such that $\rk_k(\mathbf G) \geq 2$ and set $G = \mathbf G(k)$. Let $\mathbf P < \mathbf G$ be a minimal parabolic $k$-subgroup and set $P = \mathbf P(k)$. Let $M$ be a $G$-von Neumann algebra and $E : M \to L^\infty(G/P)$ a faithful normal ucp $G$-map. The following dichotomy holds:
\begin{itemize}
\item  Either $E$ is $G$-invariant.
\item Or there exists a commutative $G$-von Neumann subalgebra $M_0 \subset M$ such that the action $G^+ \actson M_0$ is non-trivial and moreover, when writing $M_0 = L^\infty(X,\nu)$, the corresponding nonsingular action $G \actson (X,\nu)$ has the following property: for almost every $x \in X$, the stabilizer $G_x < G$ contains the $k$-points of a non-trivial $k$-split torus of $\bfG$.
\end{itemize}
\end{thm}

We note that there is no ergodicity assumption for the $G$-action on $M$ in this theorem.
The rest of this section is devoted to proving Theorem \ref{First half}. 
To this end, let us consider our almost simple $k$-group $\bfG$ with $k$-rank at least $2$, and setting $G := \bfG(k)$, let us consider a von Neumann action $G \actson M$, with a $G$-equivariant faithful normal ucp map $E: M \to L^\infty(G/P)$. Here $P = \bfP(k)$, where $\bfP$ denotes a minimal parabolic $k$-subgroup of $\bfG$.

More generally, we will freely use the notation introduced in Section \ref{AGnot}.

%Let us assume that $\bfP$ is constructed out of a maximal $k$-split torus $\bfS$ of $\bfG$ and a set of positive roots associated with this torus.

Let us assume that $E$ is not $G$-invariant. This means that $E(M)$ is not contained in the scalars. Since $G$ has rank at least $2$, the intersection of all $L^\infty(G/P_\theta)$, for $\theta \subsetneq \Delta$, is equal to the scalars inside $L^\infty(G/P)$.
Thus there exists a proper subset $\theta \subsetneq \Delta$ such that the range of $E$ is not contained in $L^\infty(G/P_\theta)$. 

For notational simplicity, we write $Q := P_\theta$, $U := V_\theta$, $H := H_\theta$, so that $Q = HU$ and $S := T_\theta$. We denote by $\oQ$ and $\oU$ the opposite parabolic group and its unipotent radical.

\begin{lem}
$E$ is not $H$-invariant on $M^{\oU}$.
\end{lem}
\begin{proof}
To prove this, denote by $B \subset M$ the maximal compact model for the $\oU$-action.
By Lemma \ref{LUdec} we note that, as a $\oU$-algebra, $L^\infty(G/P)$ is isomorphic with $L^\infty(\oU) \ovt N_0$, with $N_0 := L^\infty(Q/P)$. Moreover in this isomorphism, $L^\infty(G/Q)$ identifies with $L^\infty(\oU) \ot 1$. So, when thinking only in terms of $\oU$-algebras, we may view $E$ as a ucp map $E: M \to L^\infty(\oU) \ovt N_0$, the range of which is not contained in $L^\infty(\oU) \ovt 1$.

Restricting to compact models, and using Lemma \ref{MCM1}, we obtain a ucp map $E: B \to C_b(\overline U,N_0)$. Since $E$ is normal, we may find $b_0 \in B$, $u \in \oU$, such that $E(b_0)(u)$ is not a scalar inside $N_0$. By equivariance of $E$, the element $b _1:= \sigma_{u^{-1}}(b)$ satisfies $E(b_1)(e) \notin \C$.

Fix now a torus element $s \in S$ such that $\lim_{n \to +\infty} s^{-n} \ou s^n = e$, for all $\ou \in \oU$, and define a ucp map $E_s : M \to M$ as a point-ultraweak limit of the maps $E_s^n: x \mapsto \frac{1}{n} \sum_{k = 1}^n \sigma_{s^k}(x)$. We note that $E_s$ is not normal in general. 

{\bf Claim.} For every $b \in B$, we have $E_s(b) \in M^{\oU} =  B^{\oU}$ and $(E \circ E_s)(b) \in C_b(\oU,N_0)$ is the constant function equal to $E(b)(e)$.

For $\ou \in \oU$, the sequence $(\sigma_{s^{-n}\ou s^n}(b) - b)_{n \in \N}$, converges in norm to $0$ and hence so does the sequence of its Cesaro average. This implies that 
\[\Vert\sigma_{\ou}(E_s^n(b)) - E_s^n(b)\Vert \leq \frac{1}{n}\sum_{k = 1}^n \Vert\sigma_{s^{-k}\ou s^k}(b) - b\Vert  \to 0.\]
So a fortiori the ultraweak limit $\sigma_{\ou}(E_s(b)) - E_s(b)$ is $0$, proving the first part of the claim.\\
As explained in Lemma \ref{LUdec}, $E$ carries the $s$-action on $M$ to the diagonal action on $L^\infty(\oU) \ovt N_0$ deduced from the conjugation action $s \actson \oU$ on the one hand, and the translation action $s \actson Q/P$ on the other hand. But since $s$ centralizes $H$ and belongs to $P$, the Levi decomposition $Q = HU$ implies that $s$ acts trivially on $Q/P$. So in summary, $E \circ \sigma_s = (\sigma_s \ot \id) \circ E$.

Now, for $f \in C_b(\oU,N_0)$, the sequence $((\sigma_{s^n} \ot \id)(f))_{n \in \N}$ converges pointwise to the constant function equal to $f(e)$. By Lebesgue convergence theorem, it converges a fortiori in the ultraweak topology inherited from the embedding $C_b(\oU,N_0) \subset L^\infty(\oU) \ovt N_0$. Hence the Cesaro average of this sequence also ultraweakly converges to $f(e)$. In view of the previous paragraph, in the special case where $f = E(b)$, the Cesaro average is precisely $E \circ E_s^n(b)$. Since $E$ is normal, taking ultraweak limits, we find $E \circ E_s(b) = \lim_n E \circ E_s^n(b) = E(b)(e)$, as claimed.

We can now conclude the proof of the lemma as follows. Assume that $E$ is $H$-invariant on $M^{\oU}$. Then since $\oQ$ acts ergodically on $L^\infty(G/P)$, we conclude that in fact $E$ maps $M^{\oU}$ into $\C$. In particular, for the element $b_1 \in B$ defined above, we have $E(E_s(b_1)) \in \C$. By the claim, we have $E(E_s(b_1)) = E(b_1)(e)$. This contradicts the choice of $b_1$.
\end{proof}

Denote by $R := SU$ and $\oR := S\oU$ the solvable radical of $Q$ and $\oQ$ respectively.

\begin{lem}\label{non-invariant}
In fact, $E$ is not $H$-invariant on $M^{\oR}$.
\end{lem}
\begin{proof}
Consider the non-empty compact convex space
\[\cC:= \{\Phi: M^{\oU} \to M^{\oU} \text{ ucp map } \mid  E \circ \Phi = E|_{M^{\oU}} \text{ and } \Phi \text{ is $H$-equivariant}\},\]
with respect to the point-ultraweak topology.
Since $S$ centralizes $H$ and acts trivially on $L^\infty(G/P)^{\oU}$ (by Lemma \ref{LUdec}), we may define an action of $S$ on $\cC$ by the formula $s \cdot \Phi := \sigma_s \circ \Phi$, for all $s \in S$, $\Phi \in \cC$. This is a continuous affine action. So by amenability of $S$, it admits a fixed point $\Phi$. The range of $\Phi$ is contained in $M^{\oR}$ and $E \circ \Phi = E$ on $M^{\oU}$. Since $E$ is not $H$-invariant on $M^{\oU}$, it cannot be $H$-invariant on $\Phi(M^{\oU})$, and a fortiori on $M^{\oR}$.
\end{proof}

By Proposition \ref{embedding}, we may find a $P$-von-Neumann algebra $N$ with a $P$-invariant state $\psi$  and a normal $G$-equivariant embedding into the induced algebra $M \subset \tN := (L^\infty(G) \ovt N)^P$, such that $E$ is the restriction to $M$ of $\tE = \id \ot \psi$.

We prove that $\tN^{\oR}$ is actually contained in a nice $G$-invariant subalgebra of $\tN$. This unfortunately only holds after cutting down by a suitable projection. In the commutative setting, this projection is not relevant. This consideration already appeared in \cite{BH19}.

For this, we denote by $q \in N$ the support projection of $\psi$ restricted to $(N^R)' \cap N$. Since $R$ is normal inside $P$, $q$ is $P$-invariant. Set $p = 1 \ot q \in \tN$ and observe that $p \in \tN$ is $G$-invariant. In view of Lemma \ref{LUdec}, $\tN = \Ind_P^G(N) = \Ind_Q^G(\Ind_P^Q(N))$ may be identified with $L^\infty(\oU) \ovt N_Q$ where $N_Q = \Ind_P^Q(N) = (L^\infty(Q) \ovt N)^P$. This isomorphism maps $\tN^{\oR}$ to $1 \ot N_Q^S$ and the projection $p$ is mapped to $1 \ot (1 \ot q)$, with $1 \ot q \in N_Q$.

\begin{lem}\label{projection p}
The projection $p \in \cZ(\tN^{\oR})$ is $G$-invariant and satisfies $\tE(p) = 1$ and $p\tN^{\oR} \subset p(L^\infty(G) \ovt N^R)^P$.
\end{lem}

\begin{proof}
Using the identifications $\tN^{\oR} = 1 \ot N_Q^S$ and $p = 1 \otimes (1 \otimes q)$, to prove the lemma, it suffices to prove the following claim.

{\bf Claim.} The projection $(1 \ot q) \in N_Q$ commutes with $N_Q^S$ and $(1 \ot q)N_Q^S = (1\ot q)N_Q^R = (1 \ot q)(L^\infty(Q) \ovt N^R)^P$.

Since $R = S U$ is normal in $Q$ and contained in $P$, the subalgebra $N^R \subset N$ is $P$-invariant and $R$ acts trivially on $Q/P$. This implies that $N_Q^R = \Ind_P^Q(N)^R = \Ind_P^Q(N^R) = (L^\infty(Q) \ovt N^R)^P$ and that the ucp map $E_Q = \id \ot \psi : N_Q \to L^\infty(Q/P)$ is $R$-invariant. Moreover, $1 \ot q$ is precisely the support projection of $E_Q$ on the von Neumann subalgebra $(L^\infty(Q) \ovt ((N^R)' \cap N))^P = (N_Q^R)' \cap N_Q$. So, composing $E_Q$ with any faithful normal state on $L^\infty(Q/P)$, we get an $R$-invariant state $\psi_Q$ on $N_Q$ such that $1 \ot q$ is the support of $\psi_Q$ restricted to $(N_Q^R)' \cap N_Q$. In the GNS representation of $(N_Q,\psi_Q)$, $1 \otimes q$ is the orthogonal projection onto the closed linear span of $(N_Q)' N_Q^R\xi_{\psi_Q}$.

Let now $x \in N_Q^S$ and $g \in U$. Pick $s_n \in S$ such that $\lim_n s_ngs_n^{-1} = e$. For all $a \in (N_Q)'$ and $y \in N_Q^R$, $n \in \N$, we have
\[\Vert(\sigma_g(x) - x)ay\xi_{\psi_Q}\Vert \leq \Vert a \Vert \cdot \Vert(\sigma_g(x) - x)y\Vert_{\psi_Q} = \Vert a \Vert \cdot \Vert (\sigma_{s_ngs_n^{-1}}(x) - x)y\Vert_{\psi_Q}\]
which tends to $0$ as $n \to \infty$. This proves that $\sigma_g(x (1 \otimes q)) = \sigma_g(x)(1 \otimes q) = x (1\otimes q)$. Hence $x(1 \otimes q) \in N_Q^U \cap N_Q^S = N_Q^R$.  Since $1 \otimes q$ commutes with $N_Q^R$, we further get $x(1 \otimes q) = (1 \otimes q)x(1 \otimes q)$. The same equality for $x^*$ leads to $(1 \otimes q)x (1 \otimes q) = (1 \otimes q)x$. Hence $(1 \otimes q) \in (N_Q^S)'$ and $(1 \otimes q)N_Q^S = (1 \otimes q)N_Q^R$.

This finishes the proof of the claim and the proof of the lemma.
\end{proof}

Keep the notation $N_Q$ as in the proof of Lemma \ref{projection p} and observe that $\tN$ is naturally identified with $\Ind_Q^G(N_Q) = (L^\infty(G) \ovt N_Q)^Q$.

We now denote by $M_1 \subset M$ the $G$-invariant von Neumann subalgebra generated by $M^{\oR}$.
When viewed inside $\tN$, $M_1$ commutes with $p$, because $M^{\oR}$ commutes with $p$ and $p$ is $G$-invariant. By the previous lemma, we have $M_1p \subset (L^\infty(G) \ovt qN_Q^R)^Q$.

Since $\tE(p) = 1$ and $E$ is faithful on $M_1 \subset M$, the central support of $p \in M_1'$ inside $M_1$ is $1$. Hence the cut-down morphism $M_1 \to M_1p$ is injective. So abusing with notation, we view $M_1$ as a $G$-invariant von Neumann subalgebra of $(L^\infty(G) \ovt N_1)^Q$, with $N_1 = qN_Q^R$. 
We then define the commutative $G$-invariant von Neumann subalgebra $M_0 := \cZ(M_1) \subset M$.

\begin{lem}\label{lem-non-trivial}
The action of $G^+ \curvearrowright M_0$ is non-trivial and hence $M_0$ is non-trivial.
\end{lem}
\begin{proof}
By Lemma \ref{LUdec}, $(L^\infty(G) \ovt N_1)^Q$ is identified with $L^\infty(\oU) \ovt N_1$, and the $\oR = S\oU$-action is explicit in this description. Denote by $A_1 \subset M_1$ the maximal compact model for the $G$-action on $M_1$. It is in particular contained in the maximal compact model for the $\oU$-action, and thus, it is contained inside $C_b(\oU,N_1)$.

Denote by $N_2 \subset N_1$ the von Neumann subalgebra generated by all the values of the functions inside $A_1 \subset C_b(\oU,N_1)$. Then we have $A_1 \subset C_b(\oU,N_2) \subset L^\infty(\oU) \ovt N_2$, so that $M_1 \subset L^\infty(\oU) \ovt N_2$.

{\bf Claim 1.} $M_1$ contains $1 \ovt N_2$.

Take $f \in A_1$ and $\ou \in \oU$. It suffices to check that $1 \ot f(\ou) \in M_1$. Replacing $f$ by $\sigma_{\ou^{-1}}(f)$, we may assume that $\ou = e$. Then choose a sequence $(s_n)_{n \in \N}$ in $S$ such that $\lim_n s_n^{-1}\ov s_n = e$ for all $\overline v \in \oU$. Then the sequence $\sigma_{s_n}(f) \in A_1$ converges pointwise to $1 \ot f(e) \in C_b(\oU,N_2)$. So a fortiori, it converges to $1 \ot f(e)$ ultraweakly, which implies that $1 \ot f(e) \in M_1$, as claimed.

So thanks to the claim we have inclusions
\[1 \ot N_2 \subset M_1 \subset L^\infty(\oU) \ovt N_2.\]
In particular, $1 \ot \cZ(N_2) \subset M_0 \subset L^\infty(\oU) \ovt \cZ(N_2)$.

{\bf Claim 2.} $M_1 \neq 1 \ot N_2$. 

Indeed, otherwise $\oU$ would act trivially on $M_1$. By Tits'\ simplicity theorem \cite{Ti64}, this would imply that $G^+$ itself acts trivially on $M_1$. So in this case $E$ would be $G^+$-invariant on $M_1$. Since $G^+$ acts ergodically on $G/P$, we would conclude that $E$ is in fact $G$-invariant on $M_1$, contradicting Lemma \ref{non-invariant}. This proves Claim 2.

To conclude the lemma, it suffices to prove that $\oU$ acts non-trivially on $M_0$. Let us assume the contrary
and derive a contradiction. We could do this by combining a result of Ge-Kadison from \cite{GK95} and a direct integral argument as it was done in \cite{BH19}; we will provide a different argument, not involving direct integrals.

If $\oU$ acts trivially on $M_0$, then $M_0 = 1 \ot \cZ(N_2)$. Denote by $\Phi : N_2 \to \cZ(N_2)$ a proper normal conditional expectation. Here proper means that $\Phi(x)$ belongs to the ultraweak closure of $\{uxu^* \mid u \in \cU(N_2)\}$ for all $x \in N_2$. Such a conditional expectation exists by \cite[Theorem C]{GK95}. 
Because $\Phi$ is proper, the map $\id \ot \Phi: L^\infty(\oU) \ovt N_2 \to L^\infty(\oU) \ovt \cZ(N_2)$ maps $M_1$ into $M_1 \cap (L^\infty(\oU) \ovt \cZ(N_2)) = M_0 = 1 \ot \cZ(N_2)$. 

Take a normal faithful state $\psi$ on $\cZ(N_2)$, and denote by $\phi := \psi \circ \Phi$. Then $\phi$ is faithful on $\cZ(N_2)$. Moreover, $(\id \ot \phi)(M_1) \subset (\id \ot \psi)(1 \ot \cZ(N_2)) = \C$.

{\bf Claim 3.} The set $\Lambda := \{a\phi b \mid a,b \in N_2\}$ is norm dense in $(N_2)_*$.

The annihilator of $\Lambda$ in $N_2$ is a weak-* closed two sided ideal of $N_2$, hence of the form $zN_2$ for some $z \in \cZ(N_2)$. Since $\phi$ is faithful on $\cZ(N_2)$ we must have $z = 0$. So the claim follows from Hahn-Banach theorem.

Since $M_1$ contains $1 \ot N_2$, we find that $(\id \ot \psi')(M_1) = \C$, for every $\psi' \in \Lambda$. Hence this also holds for every $\psi' \in (N_2)_*$. It follows from \cite[Theorem B]{GK95} that $M_1 = 1 \ot N_2$. This contradicts claim 2; the proof of the lemma is complete.
\end{proof}

Denote by $A_0 \subset M_0$ the maximal compact model for the $G$-action. and denote by $X$ the Gelfand spectrum of $A_0$, so that $A_0 \simeq C(X)$. Then $M_0 = L^\infty(X,\nu)$ for some Borel measure $\nu$ on $X$\footnote{We can require that $\nu$ is a probability measure if we replace $M_0$ by a $G$-invariant separable subalgebra on which the $G^+$-action is still non-trivial.}.

\begin{lem}
Every $x \in X$ is fixed by a conjugate of the torus $S$.
\end{lem}
\begin{proof}
We want to show that every extremal state on $A_0$ is fixed by some conjugate of $S$.
We note that $A_0$ is naturally embedded in the maximal compact model $B_0$ of the $G$-von Neumann algebra $(L^\infty(G) \ovt N_1)^Q$. By Lemma \ref{C*1}, we know that every extremal state on $A_0$ is the restriction to $A_0$ of an approximately extremal state on $B_0$. Since the embedding $A_0 \subset B_0$ is $G$-equivariant, all we have to do is to check that every approximately extremal state on $B_0$ is fixed by some conjugate of $S$.

By Lemma \ref{MCM2}, we know that $B_0$ is equal to $C_b(G,C)^Q$, where $C \subset N_1$ is the maximal compact model for the $Q$-action. Since the torus $S$ is contained in the kernel of the action $Q \actson C$, the result follows from Lemma \ref{C*2}.
\end{proof}

So $M_0 = L^\infty(X,\nu)$ satisfies all the desired properties, which ends the proof of Theorem \ref{First half}.

%%%%
\section{Algebraic factors}
\label{Gauss}

Whenever $X$ is a standard Borel space, we denote by $\Prob(X)$ the space of all Borel probability measures on $X$. A standard Borel space $(X, \nu)$ endowed with a Borel probability measure is called a standard probability space. All probability spaces we consider are assumed to be standard. Recall that for any lcsc group $G$, a Borel probability measure $\mu \in \Prob(G)$ is said to be {\em admissible} if the following three conditions are satisfied:
\begin{enumerate}
\item $\mu$ is absolutely continuous with respect to the Haar measure;
\item $\supp(\mu)$ generates $G$ as a semigroup;
\item $\supp(\mu)$ contains a neighborhood of the identity element.
\end{enumerate}

In this section, we let $k$ be a local field and ${\bf G}$ an algebraic $k$-group which is assumed to be connected and absolutely almost simple (in particular, we make no rank assumption in this section). 
We will use the notation from Section \ref{AGnot}.
We consider the central isogenies $\tilde{\pi}:\tilde{\bf G}\to {\bf G}$ and $\bar{\pi}:{\bf G}\to \bar{\bf G}$, where $\tilde{\bf G}$ and
$\bar{\bf G}$ denote the simply connected cover and adjoint quotient of ${\bf G}$ correspondingly,
and let $\pi=\bar{\pi}\circ\tilde{\pi}:\tilde{\bf G}\to \bar{\bf G}$.

Our goal is to prove the following theorem.

\begin{thm} \label{thm:gauss-stat}
Assume ${\bf G}$ is absolutely almost simple.
Let $\mu$ be an admissible probability measure on $G$ and let $(X,\nu)$ be an ergodic $(G, \mu)$-stationary space.
If for a.e.\! $x\in X$ the stabilizer $G_x$ contains the $k$-points of a non-trivial split torus in ${\bf G}$, then
exactly one of the following two options holds.
\begin{itemize}
\item
The measure $\nu$ is $G$-invariant and there exists
an almost everywhere defined measurable $G$-map $G/G^+ \to X$.
\item
The measure $\nu$ is not $G$-invariant
and there exists a proper parabolic $k$-subgroup ${\bf Q}<{\bf G}$
and an almost everywhere defined measurable $G$-map $X\to ({\bf G}/{\bf Q})(k)$.
\end{itemize}
\end{thm}

Our proof of Theorem~\ref{thm:gauss-stat}
relies on some results which might be of independent interest
in the theory of algebraic groups over local fields.
It is a general fact that non-discrete closed subgroups of $G$ ``tend to be algebraic".
In Proposition~\ref{prop:ZdenseG+} we make this vague statement more precise.
In order to prove this over a local field of arbitrary characteristic
we use the above assumption that the group contains the $k$-points of a split torus,
but we comment here that this assumption could be dramatically relaxed in most situations.
However, we choose not to elaborate further on this point here.
Proposition~\ref{prop:ZdenseG+} implies 
Theorem~\ref{thm:gauss-general}, which is a non-stationary version of Theorem~\ref{thm:gauss-stat}.
Proposition~\ref{prop:cocompact} describes cocompact algebraic subgroups of $G$
and Proposition~\ref{prop:XtoV} uses the latter to describe stationary measures on $G$-algebraic varieties.
In turn, Proposition~\ref{prop:XtoV} is used to deduce Theorem~\ref{thm:gauss-stat}
from Theorem~\ref{thm:gauss-general}.
We note that the proofs below are sometimes intricate, due to phenomena of positive characteristic.

For every closed subfield of finite index, $k_0<k$,
we may perform restriction of scalars from $k$ to $k_0$ and obtain a $k_0$-algebraic group, denoted $R_{k_0}({\bf G})$,
such that $R_{k_0}({\bf G})(k_0)$ is naturally isomorphic to $G$.
Therefore, any given subgroup $L<G$ may be viewed as a subgroup of $R_{k_0}({\bf G})(k_0)$ 
and we may consider the Zariski closure of $L$ in $R_{k_0}(\bf G)$.
We will use freely the common abuse of notation, and call the Zariski closure of $L$ in $R_{k_0}(\bf G)$ its {\em $k_0$-Zariski closure in ${\bf G}$}.
In particular, we will say that $L$ is $k_0$-Zariski dense in ${\bf G}$ if it is Zariski dense in $R_{k_0}(\bf G)$.

%%%%%
%%%%%
%%%%%
%%%%%
%%%%%
%%%%%

In order to formulate properly Proposition~\ref{prop:ZdenseG+} we will need the following lemma
that, in particular, assigns to certain subgroups $L<G$ subfields $k_L<k$.

\begin{lem} \label{defn:k_L}
Assume $k$ is a non-Archimedean local field and consider a closed subgroup $L<G$ which is Zariski dense in ${\bf G}$.
Then there exists a closed subfield of finite index $k_L<k$, an absolutely simple adjoint algebraic $k_L$-group $\bar{\bf H}$
and a $k$-isogeny $\bar{\phi}:\bar{\bf H}_k\to \bar{\bf G}$ which has nowhere vanishing derivatives such that the closure of the derived subgroup 
${\tilde\pi}^{-1}(L\cap U)'$ is the image under $\tilde{\phi}:\tilde{\bf H}_k\to \tilde{\bf G}$ of an open subgroup of 
$\tilde{\bf H}(k_L)<\tilde{\bf H}(k)=\tilde{\bf H}_k(k)$.
\end{lem}

Here $\bar{\bf H}_k$ and $\tilde{\bf H}_k$ stand for the $k$-groups obtained from $\bar{\bf H}$ and $\tilde{\bf H}$ by extending the scalars from $k_L$ to $k$, where $\tilde{\bf H}$ is the simply connected cover of ${\bf H}$, 
$\pi':\tilde{\bf H}\to\bar{\bf H}$ the corresponding central $k_L$-isogeny
and $\tilde{\phi}$ is the pull back of $\phi$ along $\pi'$.

\begin{proof}
Using van Dantzig's Theorem, we choose a compact open subgroup $U$ of the totally disconnected locally compact group $G={\bf G}(k)$
and note that $L\cap U<G$ is Zariski dense by \cite[Lemma 7.5]{BF18},
as it is commensurated by $L$.
We are thus able to apply \cite[Theorem 0.2(a)]{Pi97} to the Zariski dense compact subgroup $\bar{\pi}(L\cap U)<\bar{\bf G}(k)$ 
and deduce that there exists a closed subfield of finite index $k_L<k$, an absolutely simple adjoint algebraic $k_L$-group $\bar{\bf H}$
and a $k$-isogeny $\bar{\phi}:\bar{\bf H}_k\to \bar{\bf G}$ which has nowhere vanishing derivatives such that the closure of the derived subgroup 
${\tilde\pi}^{-1}(L\cap U)'$ is the image under $\tilde{\phi}:\tilde{\bf H}_k\to \tilde{\bf G}$ of an open subgroup of 
$\tilde{\bf H}(k_L)<\tilde{\bf H}(k)=\tilde{\bf H}_k(k)$.
\end{proof}

The subfield $k_L<k$ is independent of the choice of the compact open subgroup $U<G$ by \cite[Theorem 0.2(b)]{Pi97} and commensurability. 
We thus get a well defined $G$-equivariant map $L\mapsto k_L$ from the space of closed Zariski dense subgroups
of $G$ to the space of closed subfields of finite index of $k$.

\begin{prop} \label{prop:ZdenseG+}
Assume ${\bf G}$ is absolutely almost simple.
Let $L<G$ be a proper closed subgroup which contains the $k$-points of a non-trivial split torus in ${\bf G}$
and does not contain $G^+$.
Then there exists a closed subfield of finite index $k_0<k$ such that $L$ is not $k_0$-Zariski dense. 

We obtain the following more precise statement.

\begin{itemize}
\item If $L$ is not Zariski dense in ${\bf G}$, then we can clearly take $k_0=k$.

\item If $L$ is Zariski dense in ${\bf G}$, then in case $k$ is Archimedean, we can take $k_0=\mathbb{R}$.
Otherwise, we can take $k_0=k_L$, where $k_L$ is the subfield given in Lemma~\ref{defn:k_L}, unless the characteristic of $k$ is 2 or 3,
in which case we can take either $k_0=k_L$ or $k_0=F(k)$, where $F:k\to k$ is the Frobenius endomorphism.
\end{itemize}
\end{prop}

\begin{proof}
We assume as we may that $L$ is Zariski dense in ${\bf G}$.
Note that by \cite[Theorem I.2.3.1(a) and Proposition I.1.5.5]{Ma91} we have that $G^+=\tilde{\pi}(\tilde{\bf G}(k))$
and we recall that this is a closed, cocompact subgroup.
We set $\tilde{G}=\tilde{\bf G}(k)$, $\bar{G}=\bar{\bf G}(k)$, $\bar{G}^+=\pi(\tilde{G})$, $\tilde{L}=\tilde{\pi}^{-1}(L)$ and $\bar{L}=\bar{\pi}(L)$,
thus $\tilde{L}$ and $\bar{L}$ are Zariski dense in $\tilde{\bf G}$ and $\bar{\bf G}$ correspondingly.
Since the restriction of scalar functor applied to a central isogeny is a central isogeny, 
it is enough to show that $\bar{L}$ is not $k_0$-Zariski dense in $\bar{\bf G}$.
We proceed to show that.

In case $k$ is Archimedean the proof is quite standard:
$\bar{L}$ is not $\mathbb{R}$-Zariski dense in $\bar{\bf G}$,
as it is contained in the stabilizer of the Lie-algebra of $L$
in the real adjoint representation of $G$.
Here we consider $G$ as a real Lie group and use the fact that 
$G^+$ is the identity connected component of $G$ as such, thus the Lie algebra of $L$ is 
a non-trivial sub Lie-algebra of the Lie-algebra of $G$, and in particular it is not $G$-invariant,
as the adjoint representation is irreducible.

From now on we will assume as we may that $k$ is non-Archimedean
and consider the field $k_L$ given in Lemma~\ref{defn:k_L}.
We will show that if $k_L\lneq k$ is a proper subfield then
$\bar{L}$ is not $k_L$-Zariski dense in $\bar{\bf G}$ 
and if $k_L=k$ then necessarily the characteristic of $k$ is 2 or 3 and 
$\bar{L}$ is not $F(k)$-Zariski dense in $\bar{\bf G}$.
This will finish the proof.

We fix a compact open subgroup $U<G$
and consider, along with the field $k_L$, the absolutely simple adjoint algebraic $k_L$-group $\bar{\bf H}$
and the $k$-isogeny $\bar{\phi}:\bar{\bf H}_k\to \bar{\bf G}$ which has nowhere vanishing derivatives such that 
${\tilde\pi}^{-1}(L\cap U)'$ is the image under $\tilde{\phi}:\tilde{\bf H}_k\to \tilde{\bf G}$ of an open subgroup of 
$\tilde{\bf H}(k_L)<\tilde{\bf H}(k)=\tilde{\bf H}_k(k)$.
We set $M=L\cap U$, $\bar{M}=\bar{\pi}(M)<\bar{L}$ and $\tilde{M}=\tilde{\pi}^{-1}(M)<\tilde{L}$.
Note that ${\tilde\pi}^{-1}(L\cap U)'=\tilde{M}'$ is the image under $\tilde{\phi}:\tilde{\bf H}_k\to \tilde{\bf G}$ of an open subgroup of 
$\tilde{\bf H}(k_L)$.

We now claim that if $k_L\lneq k$ is a proper subfield then $\bar{L}$ is not $k_L$-Zariski dense in $\bar{\bf G}$. By contradiction, assume that $k_L\lneq k$ and $\bar{L}$ is $k_L$-Zariski dense in $\bar{\bf G}$.
Note that $\bar{M}$ is a compact open subgroup of $\bar{L}$,
thus $\bar{L}$ commensurates $\bar{M}$.
By the $k_L$-Zariski density assumption on $\bar{L}$, we get by \cite[Lemma 7.5]{BF18} that $\bar{M}$ is $k_L$-Zariski dense in $\bar{\bf G}$.
Thus by the $k_L$-simplicity of $\bar{\bf G}$ we deduce that $\bar{M}'$ is $k_L$-Zariski dense in $\bar{\bf G}$ as well.
However, since $k_L\lneq k$ we get that 
$\tilde{\phi}^{-1}(\tilde{M}')$ is not $k_L$-Zariski dense in $\tilde{\bf H}_k$,
as the $k_L$-points of its $k_L$-Zariski closure are contained in $\tilde{\bf H}(k_L)$
which is a proper subgroup of $\tilde{\bf H}(k)=\tilde{\bf H}_k(k)$
and it follows that $\bar{M}'$ is not $k_L$-Zariski dense in $\bar{\bf G}$.
This gives the desired contradiction.

From now on we assume as 
we may that $k_L=k$ and argue to show that necessarily the characteristic of $k$ is 2 or 3 and 
$\bar{L}$ is not $F(k)$-Zariski dense in $\bar{\bf G}$.
We have that $\tilde{\phi}^{-1}(\tilde{L})$ is open in $\tilde{H}=\tilde{\bf H}(k)$, as it contains $\tilde{\phi}^{-1}(\tilde{M}')$,
and it is non-compact, as it contains the $k$-points of the $\tilde{\phi}\circ\tilde{\pi}$-preimage of a split torus of ${\bf G}$.
It follows by Howe-Moore Theorem that $\tilde{\phi}^{-1}(\tilde{L})=\tilde{H}$ (see \cite[Theorem 5.1]{HM77}).
We conclude that $\bar{\phi}$ is not an isomorphism.
Indeed, if it was then also $\tilde{\phi}$ was an isomorphism 
and we would get 
$G^+=\tilde{\pi}(\tilde{G})=\tilde{\pi}\circ\tilde{\phi}(\tilde{H})=\tilde{\pi}(\tilde{L})<L$,
which contradicts our assumption that $L$ does not contain $G^+$.
We are thus in a very special situation where $\bar{\phi}$ is a 
$k$-isogeny which is not an isomorphism which has nowhere vanishing derivatives.
Such non-standard isogenies are discussed in \cite[\S1]{Pi97}.
They appear only in a few special cases: either the characteristic of $k$ is 3 and the root systems of $\bar{\bf G}$
and $\bar{\bf H}$ are of type $G_2$ or 
the characteristic of $k$ is 2 and the root systems of $\bar{\bf G}$
and $\bar{\bf H}$ are of types $B_n$, $C_n$ or $F_4$ and $C_n$, $B_n$ or $F_4$ correspondingly.
In particular, if the characteristic of $k$ is not 2 or 3 then such non-standard isogenies
do not exist.
We thus have that the characteristic of $k$ is 2 or 3 and the types of $\bar{\bf G}$ and $\bar{\bf H}$
are as above and we are left to show that $\bar{L}$ is not $F(k)$-Zariski dense.

By \cite[Propositions 1.6 and 1.7]{Pi97} 
there exists an absolutely simple adjoint $k$-group $\bar{\bf H}^\sharp$ and isogenies
\[ \psi:\bar{\bf H}\to \bar{\bf H}^\sharp, \quad \psi^\sharp:\bar{\bf H}^\sharp\to \bar{\bf H}_F
\quad \mbox{and} \quad
\chi: \bar{\bf H}^\sharp \to \bar{\bf G} \]
such that $\bar{\phi}=\chi\circ \psi$ and $\psi^\sharp\circ \psi=F^*:\bar{\bf H}\to \bar{\bf H}_F$,
where $\bar{\bf H}_F$ is the $k$-group obtained from the $k$-group $\bar{\bf H}$ by extending scalars via $F:k\to k$
and $F^*$ is the Frobenius isogeny.
Setting $\theta=\psi^\sharp\circ\chi^{-1}$ we have 
$F^*=\psi^\sharp\circ \psi=\psi^\sharp\circ \chi^{-1}\circ \chi\circ \psi=\theta\circ \bar{\phi}$
and thus
obtain the following diagram
which presents some of the groups and isogenies relevant to our 
discussion.
%\[
%\begin{tikzcd}
%\tilde{\bf H} \arrow{d}{\tilde{\phi}} \arrow{r}{\pi'}  & \bar{\bf H} \arrow{d}{\bar{\phi}} \arrow[bend left=60]{dd}{F^*} \\
%\tilde{\bf G} \arrow{r}{\pi}  & \bar{\bf G} \arrow{d}{\theta} \\
%& \bar{\bf H}_F
%\end{tikzcd}
%\]
\begin{equation} \label{eq:isogenies}
\begin{tikzcd}
\tilde{\bf H} \arrow{r}{\tilde{\phi}} \arrow{d}{\pi'}  & \tilde{\bf G} \arrow{d}{\pi} \\
\bar{\bf H} \arrow{r}{\bar{\phi}} \arrow[bend right=45]{rr}{F^*} & \bar{\bf G} \arrow{r}{\theta} & \bar{\bf H}_F
\end{tikzcd}
\end{equation}

We note that, by the construction of $\psi$ in \cite[Propositions 1.6]{Pi97},
we have that $\psi(\bar{H})$ is a proper subgroup of $\bar{H}^\sharp=\bar{\bf H}^\sharp(k)$,
thus $\bar{\phi}(\bar{H})=\chi\circ\psi(\bar{H})$ is a proper subgroup of $\bar{G}=\chi(\bar{H}^\sharp)$.
We claim further that  $\bar{\phi}(\bar{H})$ is not $F(k)$-Zariski dense in $\bar{G}$.
Applying $\theta$, it is enough to show that $F^*(\bar{H})=\theta\circ \bar{\phi}(\bar{H})$
is not $F(k)$-Zariski dense in $\theta(\bar{G})$.
Since $\theta$ is injective, it is thus enough to show that the group $F^*(\bar{H})$
equals the group of $F(k)$-points of its $F(k)$-Zariski closure in $\bar{\bf H}_F$.
Therefore, we need to understand $R_{F(k)}(\bar{\bf H}_F)$ and the image of $R_{F(k)}(F^*)$ in it.
We construct the extension of scalars of $\bar{\bf H}$ via $F:k\to k$ in two stages:
we define $\bar{\bf H}_{F(k)}$ to be the $F(k)$-group obtained by extending the scalars 
of $\bar{\bf H}$ via $F:k\to F(k)$ and we view $\bar{\bf H}_F$ as the 
$k$-group obtained by extending the scalars 
of $\bar{\bf H}_{F(k)}$ via the inclusion $F(k)\to k$.
The adjunction of restriction and extension of scalars functors provides a canonical inclusion,
called the unit map, of $F(k)$-groups $\bar{\bf H}_{F(k)} \to R_{F(k)}(\bar{\bf H}_F)$
and we identify $\bar{\bf H}_{F(k)}$ with its image in $R_{F(k)}(\bar{\bf H}_F)$.
By the definition of $F^*$, we have that the image of $R_{F(k)}(\bar{\bf H})$
under $R_{F(k)}(F^*)$ is $\bar{\bf H}_{F(k)}$,
the map being the application of $F$ on the $k$-coordinates.
In particular, on the level of $F(k)$-points, $R_{F(k)}(F^*)$ provides an isomorphism
$R_{F(k)}(\bar{\bf H})(F(k))\simeq \bar{\bf H}(k)\simeq \bar{\bf H}_{F(k)}(F(k))$.
We conclude that the $F(k)$-Zariski closure of $F^*(\bar{H})$ is $\bar{\bf H}_{F(k)}$
and its group of $F(k)$-points is again $F^*(\bar{H})$.
This finishes the proof of the claim that $\bar{\phi}(\bar{H})$ is not $F(k)$-Zariski dense in $\bar{G}$.

We denote $\bar{N}=\pi\circ\tilde{\phi}(\tilde{H})=\bar{\phi}\circ \pi(\tilde{H})$
and conclude that it is not $F(k)$-Zariski dense in $\bar{\bf G}$, as $\bar{N}<\bar{\phi}(\bar{H})$.
Since $\tilde{\phi}^{-1}(\tilde{L})=\tilde{H}$,
we have that $\tilde{\phi}(\tilde{H})<\tilde{L} < \pi^{-1}(\bar{L})$, thus $\bar{N}<\bar{L}$.
We claim that $\bar{N}$ is normal in $\bar{L}$.
Denoting $\bar{H}=\bar{\bf H}(k)$ and $\bar{H}^+=\pi'(\tilde{H})<\bar{H}$,
we have that $\bar{H}^+$ is normal in $\bar{H}$.
Thus $\bar{N}=\bar{\phi}\circ\pi'(\tilde{H})=\bar{\phi}(\bar{H}^+)$
is normal in $\bar{\phi}(H)$.
Therefore, we will be done by showing that $\bar{L}<\bar{\phi}(\bar{H})$.
To show this, we fix $\bar{g}\in \bar{L}$ and argue to show that $\bar{g}\in \bar{\phi}(\bar{H})$.
We will denote by $\iota=\inn_{\bar{g}}$ 
the corresponding conjugation automorphism of $\bar{\bf G}$
and we will denote by $\tilde{\iota}$ the corresponding $k$-automorphism of $\tilde{\bf G}$.
Since $(\tilde{\iota}\circ \tilde{\phi})^{-1}(\tilde{\iota}(\tilde{M}'))=\tilde{\phi}^{-1}(\tilde{M}')$ is open in $\tilde{\bf H}(k)$,
we get that the triple $(k,\bar{\bf H},\iota\circ\bar{\phi})$ satisfies \cite[Theorem 0.2(a)]{Pi97} 
for the compact group $\iota(\bar{M})$ 
as well as for its finite index subgroup $\iota(\bar{M}) \cap \bar{M}$
(we remind the reader that $\bar{M}$ is commensurated in $\bar{L}$).
%But, since the triple $(k,\bar{\bf H},\bar{\phi})$ forms a solution for the group $\bar{M}$,
%it also forms such a solution for its finite index subgroup $\iota(\bar{M}) \cap \bar{M}$.
Thus, 
we obtain by the uniqueness statement of \cite[Theorem 0.2(b)]{Pi97}
that there exists a $k$-automorphism $\alpha$ of $\bar{\bf H}$ such that $\iota\circ\bar{\phi}=\bar{\phi}\circ \alpha$.
However, $\bar{\bf H}$ is an adjoint group and its Dynkin diagram
is of type $B_n, C_n, F_4$ or $G_2$, thus have no non-trivial automorphism. 
Therefore, every automorphism of $\bar{\bf H}$ is inner and we conclude that $\alpha=\inn_{\bar{h}}$ for some $\bar{h}\in \bar{H}$.
We thus get 
\[ \inn_{\bar{g}}\circ\bar{\phi}=\iota\circ\bar{\phi}=\bar{\phi}\circ \alpha=\bar{\phi}\circ \inn_{\bar{h}}=\inn_{\bar{\phi}(\bar{h})}\circ\bar{\phi}\]
and therefore $\inn_{\bar{g}}=\inn_{\bar{\phi}(\bar{h})}$ on $\bar{\bf G}$,
as $\bar{\phi}$ is dominant.
Since $\bar{\bf G}$ is adjoint we conclude that indeed $\bar{g}=\bar{\phi}(\bar{h})\in \bar{\phi}(\bar{H})$.
This finishes the proof of the claim that $\bar{N}$ is normal in $\bar{L}$.

Since $\bar{N}$ is not $F(k)$-Zariski dense in $\bar{\bf G}$ and it is normal in $\bar{L}$,
we conclude by $F(k)$-simplicity that $\bar{L}$ is not $F(k)$-Zariski dense in $\bar{\bf G}$.
This finishes the proof.
\end{proof}

\begin{thm} \label{thm:gauss-general}
Assume ${\bf G}$ is absolutely almost simple.
Let $X$ be an ergodic $G$-space.
If for a.e.\! $x\in X$ the stabilizer $G_x$ 
contains the $k$-points of a non-trivial split torus in ${\bf G}$,
then
exactly one of the following two options holds.
\begin{itemize}
\item
There exists 
an almost everywhere defined measurable $G$-map $G/G^+ \to X$.
\item
There exists an infinite closed subfield $k_0<k$, a proper $k_0$-algebraic subgroup ${\bf H}<R_{k_0}({\bf G})$
and an almost everywhere defined measurable $G$-map $X\to (R_{k_0}({\bf G})/{\bf H})(k_0)$.
\end{itemize}
\end{thm}

\begin{proof}
If $G^+$ is in the kernel of the $G$-action on $X$ 
then $X$ is an ergodic $G/G^+$-space
and by the compactness of $G/G^+$ the action is in fact transitive.
In that case the first bullet clearly holds.
We thus assume from now on that $G^+$ is not contained in the kernel of the $G$-action on $X$ 
and argue to show that the second bullet holds.

We let $\Sub(G)$ be the space of closed subgroup of $G$ endowed with the Chabauty topology 
and consider the stabilizer map 
\[ X \to \Sub(G), \quad x \mapsto G_x, \]
taking a point in $x$ to its stabilizer in $G$.
This map is measurable by \cite[Ch.\ II, Prop.\ 2.3]{AM66}.
%where AM= L. AUSLANDER and C. C. MOORE, Unitary Representations of Solvable Lie Groups, Memoirs of the A. M. S. 62 (1966)
%but I don't have access.
By ergodicity we deduce that on a full measure subset of $X$, $G^+$ is not contained in $G_x$.
Without loss of the generality, we will assume that this is the case for every $x\in X$.
Similarly, we assume that for every $x$, $G_x$ is a closed subgroup of $G$ which  
contains the $k$-points of a non-trivial split torus in ${\bf G}$.
By Proposition~\ref{prop:ZdenseG+} we get that for every $x\in X$ 
there exists a closed subfield of finite index $k_x<k$ such that $G_x$ is not $k_x$-Zariski dense.

Next, we claim that there exists a single closed subfield of finite index $k_0<k$ such that 
for almost every $x\in X$, $G_x$ is not $k_0$-Zariski dense in ${\bf G}$.
If $k$ is Archimedean this follows by setting $k_0=\mathbb{R}$.
We now assume that $k$ is non-Archimedean and prove this claim.
We use \cite[Proposition 1.11]{Pi97} to fix a $k$-representation $\rho$ of ${\bf G}$ 
which is a non-constant absolutely irreducible subquotient of the adjoint representation.
We consider the map $\Sub(G)\to \Sub(k)$ taking a closed subgroup of $G$ to the closed subfield of $k$ generated by 
all traces of all $\rho(g)$, where $g$ ranges over all elements of all the compact open subgroups of $L$.
As in the proof of \cite[Proposition 5.4]{GL17}, this map is Borel measurable.
It is clearly $G$-invariant map.
We conclude, by ergodicity, that the composed map $X\to \Sub(G)\to \Sub(k)$ is essentially constant.
We assume as we may that this map is actually constant and denote its value by $k_1$.
We note that by \cite[Lemma 7.5]{BF18}, for every $x\in X$, every compact open subgroup of $G_x$ is Zariski dense in ${\bf G}$.
From \cite[Proposition 0.6(a)]{Pi97} we conclude that if the characteristic of $k$ is not 2 or 3 then 
for every $x\in X$, $k_1=k_{G_x}$ and we conclude from Proposition~\ref{prop:ZdenseG+} 
that in this case $G_x$ is not $k_1$-Zariski dense, which proves the claim upon setting $k_0=k_1$.
We assume now that the charactersitic of $k$ is 2 or 3.
By Proposition~\ref{prop:ZdenseG+} we have that for every $x\in X$, $G_x$ is not $F(k_{G_x})$-Zariski dense. 
By \cite[Proposition 0.6(a)]{Pi97} we have that for every $x\in X$, 
$F(k_{G_x})<k_1<k_{G_x}$.
From $k_1<k_{G_x}$ we get that for every $x\in X$, $G_x$ is not $F(k_1)$-Zariski dense
and from $F(k_{G_x})<k_1$ we get that $k_1$, thus also $F(k_1)$, is a closed subfield of finite index in $k$.
The claim now follows by setting $k_0=F(k_1)$.

From now on we will identify $G$ with the group of $k_0$-points in the restriction of scalars of ${\bf G}$
from $k$ to $k_0$, namely $R_{k_0}({\bf G})(k_0)$.
We consider the $k_0$-algebraic group $R_{k_0}({\bf G})$ and the corresponding Zariski closure map,
\[ z:\Sub(G) \to \Sub(G), \quad H \mapsto \bar{H}^Z(k_0), \]
taking a closed subgroup of $G$ to the group of $k_0$-points of its Zariski closure in $R_{k_0}({\bf G})$.
Noting that $R_{k_0}({\bf G})$ is $k_0$-almost simple,
we have that the only normal proper $k_0$-subgroups of $R_{k_0}({\bf G})$ are central, (hence finite).
We conclude that the map $X\to \Sub(G)$, $x\mapsto z(G_x)$ is not essentially constant.

Next, we fix a faithful $k_0$-linear representation $R_{k_0}({\bf G})\to \SL_n$ and consider $R_{k_0}({\bf G})$ as a subvariety of the space of matrices,
which we identify with the vector space $k_0^{n^2}$.
We thus identify $\Sub(G)$ as a closed subspace of $\Cl(k_0^{n^2})$,
the space of closed subsets of $k_0^{n^2}$. 
We set $A=k_0[x_1,\dots,x_{n^2}]$
and consider the map taking $x\in X$ to the ideal $I_x\lhd A$ consisting of the polynomials vanishing on $G_x<k_0^{n^2}$.
We conclude that this map is not essentially constant.

We denote by $A_i<A$ the linear subspace consisting of polynomials of degree bounded by $i$.
We consider the ``chopping map" $\sigma_i:A\to A_i$, which takes polynomial to the sum of its monomials of degree bounded by $i$.
We let $\Gr(A_i)=\cup_{d=0}^{\dim(A_i)} \Gr(d,A_i)$ be the full Grassmannian of $A_i$
and consider the map
\[ \psi:X \to \Gr(A_i), \quad x \mapsto \sigma_i(I_x). \]
This is a composition of the stabilizer map $X\to \Cl(k_0^{n^2})$, $x\mapsto G_x$ and the map $\Cl(X) \to \Gr(A_i)$, 
taking $F\subset k_0^{n^2}$ to the image under $\sigma_i$ of its annihilating ideal.
As noted above, the first map is measurable by \cite[Ch.\ II, Prop.\ 2.3]{AM66}.
The second map is measurable by \cite[Proposition 4.2]{GL17}.
%INVARIANT RANDOM SUBGROUPS OVER NON-Archimedean LOCAL FIELDS
We conclude that $\psi$ is measurable as well.
Picking $i$ large enough we make sure that $\psi_i$ is not essentially constant
and we set $E=A_i$ and $\psi=\psi_i$.
By ergodicity, we conclude that the image of $\psi$ is contained in $\Gr(d,E)$ for some $0< d< \dim(E)$.
We let ${\bf V}$ be the $k_0$-algebraic variety corresponding to $\Gr(d,E)$.
By \cite[Proposition 4.2]{BF18} % SUPER-RIGIDITY AND NON-LINEARITY FOR LATTICES IN PRODUCTS URI BADER AND ALEX FURMAN
we get a $G$-map $\phi:X \to (R_{k_0}({\bf G})/{\bf H})(k_0)$, for some $k_0$-algebraic subgroup ${\bf H}<R_{k_0}({\bf G})$,
such that $X\to {\bf V}(k_0)$ factors via $\phi$.
As $\psi$ is not essentially constant, we
conclude that ${\bf H}$ is a proper $k_0$-subgroup of $R_{k_0}({\bf G})$.
This finishes the proof.
\end{proof}

\begin{prop} \label{prop:cocompact}
Let ${\bf G}$ be a $k$-isotropic almost $k$-simple connected algebraic $k$-group and ${\bf H}$ a $k$-subgroup of ${\bf G}$.
Then either ${\bf H}^0$,
the identity component of ${\bf H}$, is reductive
or 
there exists a proper parabolic $k$-subgroup ${\bf Q}<{\bf G}$
such that ${\bf H}<{\bf Q}$.
If ${\bf H}$ is a proper subgroup of ${\bf G}$ and ${\bf H}(k)$ is cocompact in ${\bf G}(k)$ then 
${\bf H}^0$ is not reductive, thus the second option holds.
\end{prop}

\begin{proof}
Assume first that ${\bf H}^0$ is not reductive, thus its unipotent radical ${\bf U}$ is non-trivial
and ${\bf H}$ is contained in its normalizer.
By \cite[Corollary 3.9(ii)]{BT70}, applied over $\bar{k}$, the algebraic closure of $k$, there exists a proper parabolic subgroup 
${\bf Q}<{\bf G}$ which contains ${\bf H}$.
We need to show that ${\bf Q}$ is defined over $k$.
This follows directly from \cite[Corollary 3.9(ii)]{BT70} in case
${\bf U}$ is defined over $k$,
which is always the case if $k$ is of characteristic 0. 
For the general case, we use \cite[Corollary~18.8]{Bo91} to deduce that ${\bf G}$ is split over $k_s$,
the separable closure of $k$ in $\bar{k}$,
thus ${\bf Q}$ is defined over $k_s$.
By the fact that ${\bf H}$ is defined over $k$, we get that ${\bf U}$
is stable under all automorphisms of $\bar{k}$ which preserve $k$, $\Aut_k(\bar{k})$.
Applying \cite[Corollary 3.9(iii)]{BT70} we have that ${\bf Q}$ is invariant under $\Aut_k(\bar{k})$ as well.
Identifying $\Aut_k(\bar{k})\simeq \Aut_k(k_s)$ with the Galois group of $k_s$ over $k$,
we conclude by \cite[Theorem~AG14.4]{Bo91} that indeed, ${\bf Q}$ is defined over $k$.

Next we will assume that ${\bf H}$ is a proper subgroup of ${\bf G}$,
${\bf H}^0$ is reductive and $H={\bf H}(k)$ is cocompact in $G$.
We will argue to show a contradiction.
By the main result of \cite{Ri75} and \cite{Ha74} 
% AFFINE COSET SPACES OF REDUCTIVE ALGEBRAIC GROUPS R. W. RICHARDSON
% HOMOGENEOUS VECTOR BUNDLES  AND REDUCTIVE SUBGROUPS OF REDUCTIVE ALGEBRAIC GROUPS. By W. J. HABouSH.
we have that the variety ${\bf G}/{\bf H}$ is affine.
We note that ${\bf G}(k)$ is Zariski dense in ${\bf G}$, thus its image $B={\bf G}(k)/{\bf H}(k)$ 
is Zariski dense in ${\bf G}/{\bf H}$.
We conclude that $\Stab_{\bf G}(\bar{B}^Z)={\bf G}$ and thus ${\bf N}=\Fix_{\bf G}(\bar{B}^Z)$ is normal in ${\bf G}$.
By the $k$-almost simplicity of ${\bf G}$, we conclude that ${\bf N}$ is central in ${\bf G}$, hence finite.
Using \cite[Proposition~4.10]{BDL15} we conclude that $G$ is compact.
It follows that ${\bf G}$ is $k$-anisotropic, which gives the desired contradiction.
\end{proof}

\begin{prop} \label{prop:XtoV}
Let ${\bf G}$ be a $k$-isotropic almost $k$-simple connected algebraic $k$-group. Let
${\bf H}<{\bf G}$ be a $k$-subgroup such that $({\bf G}/{\bf H})(k)$ admits a $\mu$-stationary probability measure $\nu$
with respect to some admissible probability measure $\mu$ on $G$.
Then there exists a proper parabolic $k$-subgroup ${\bf Q}<{\bf G}$
such that ${\bf H}<{\bf Q}$.
\end{prop}

\begin{proof}
Set $G = \mathbf G(k)$ and $H = \mathbf H(k)$. Let ${\bf P}<{\bf G}$ be a minimal parabolic $k$-subgroup
and set $P={\bf P}(k)$.
By \cite[Corollary 5.2]{BS04} there is a $\mu$-stationary measure $\nu_0$ on $G/P$ such that $(G/P,\nu_0)$ 
forms the Furstenberg-Poisson boundary of $(G,\mu)$.
By \cite[Corollary 2.16]{BS04} we have a ``boundary map" $\beta:G/P\to \Prob(({\bf G}/{\bf H})(k))$
such that $\nu$ is the barycenter of $\beta_*\nu_0$.
We denote by $\xi$ the image of the base coset under $\beta$.
Thus, $\xi$ is a $P$-invariant probability measure on $({\bf G}/{\bf H})(k)$.

Next, we apply \cite[Theorem 1.1]{Sh97} 
(see also \cite[Proposition~1.9]{BDL15})
for the $k$-algebraic group ${\bf P}$ to find a normal $k$-subgroup ${\bf P}_0<{\bf P}$ such that $P/P_0$ is compact,
where $P_0={\bf P}_0(k)$, such that the measure $\xi$ is supported on the variety of ${\bf P}_0$-fixed points in ${\bf G}/{\bf H}$.
In particular, we deduce that this variety is non-empty.
Upon conjugating ${\bf H}$, we may thus assume that ${\bf P}_0<{\bf H}$.
Since $P_0<G$ is cocompact, as both $P_0<P$ and $P<G$ are, it follows that $H<G$ is cocompact.
We are now done, by Proposition~\ref{prop:cocompact}.
\end{proof}

\begin{proof}[Proof of Theorem~\ref{thm:gauss-stat}]
We consider the dichotomy provided by Theorem~\ref{thm:gauss-general}.
If $X$ is a factor of $G/G^+$ then it carries an invariant probability measure, the image of the Haar measure of $G/G^+$,
which must coincide with $\nu$ by ergodicity.
Then we get the first case of Theorem~\ref{thm:gauss-stat}.
We thus assume as we may that 
there exists an infinite closed subfield $k_0<k$, a proper algebraic $k_0$-subgroup ${\bf H}<R_{k_0}({\bf G})$
and an almost everywhere defined measurable $G$-map $X\to (R_{k_0}({\bf G})/{\bf H})(k_0)$.
We endow $(R_{k_0}({\bf G})/{\bf H})(k_0)$ with the pushforward stationary measure.
Using Proposition~\ref{prop:XtoV}
we find a proper parabolic $k_0$-subgroup ${\bf Q}_0<R_{k_0}({\bf G})$
such that ${\bf H}<{\bf Q}_0$.
We thus obtain a $k_0$-quotient map $R_{k_0}({\bf G})/{\bf H}\to R_{k_0}({\bf G})/{\bf Q}_0$
and the resulting composition 
\[ X\to (R_{k_0}({\bf G})/{\bf H})(k_0)\to (R_{k_0}({\bf G})/{\bf Q}_0)(k_0). \]
By the discussion in \cite[\S16.2.6]{Sp98}
%Here we need a reference. I found it in the following note by Brian Conrad:
%http://virtualmath1.stanford.edu/~conrad/249BW16Page/handouts/relbruhat.pdf
%Lemma 6.1.
%Possibly we can find it in his book "pseudo reductive groups"
there exists a parabolic $k$-subgroup ${\bf Q}<{\bf G}$ such that ${\bf Q}_0=R_{k_0}({\bf Q})$.
Since $(R_{k_0}({\bf G})/{\bf Q}_0)(k_0)=R_{k_0}({\bf G})(k_0)/{\bf Q}_0(k_0)={\bf G}(k)/{\bf Q}(k)=({\bf G}/{\bf Q})(k)$,
we get the desired map $X\to ({\bf G}/{\bf Q})(k)$.
The proof is complete, noting that the space $({\bf G}/{\bf Q})(k)$ admits no $G$-invariant measure,
as any such measure must be supported on ${\bf G}$-invariant points.
\end{proof}

%%%%%%%%%%

\section{Proofs of the main results}
\label{main results}

\subsection{Uniqueness of equivariant measurable maps}

In this subsection, we prove a useful fact regarding uniqueness of equivariant measurable maps. Following \cite{BFGW12}, whenever $G$ is a locally compact second countable group, $(X, \nu)$ a standard probability $G$-space and $Z$ a standard Borel $G$-space, we denote by $\Map_G(X, Z)$ the set of all equivalence classes of measurable $G$-maps $\zeta : X \to Z$.

\begin{prop}\label{prop-uniqueness}
Let $k$ be a local field. Let $\mathbf G$ be a $k$-isotropic almost $k$-simple connected algebraic $k$-group. Let $\mathbf P < \mathbf Q < \mathbf G$ be proper parabolic $k$-subgroups so that $\mathbf P < \mathbf G$ is a minimal parabolic $k$-subgroup. Endow the homogeneous space $(\bfG/\bfP)(k) = \bfG(k)/\bfP(k) $ with its unique $\bfG(k)$-invariant measure class and denote by $p : \bfG(k)/\bfP(k)\to \bfG(k)/\bfQ(k)$ the canonical $\bfG(k)$-map. 

Then for every lattice $\Gamma < \bfG(k)$, we have
$$\Map_{\bfG(k)}(\bfG(k)/\bfP(k) , \Prob(\bfG(k)/\bfQ(k) )) = \{x \mapsto \delta_{p(x)}\} = \Map_\Gamma(\bfG(k)/\bfP(k) , \Prob(\bfG(k)/\bfQ(k))).$$
\end{prop}

\begin{proof}
Denote by $p_* : \Prob(\bfG(k)/\bfP(k)) \to \Prob(\bfG(k)/\bfQ(k))$ the pushforward $\bfG(k)$-map. 

Firstly, assume that $\rk_k(\bfG) = 1$. Then $\Gamma \curvearrowright \bfG(k)/\bfP(k)$ is a convergence action in the sense of \cite[Section 3]{BF14} and \cite[Theorem 3.2]{BF14} implies that $ \bfG(k)/\bfP(k) \to  \Prob(\bfG(k)/\bfQ(k)) : x\mapsto \delta_{p(x)}$  is the essentially unique measurable $\Gamma$-map (resp.\ $\bfG(k)$-map).

Secondly, assume that $\rk_k(\bfG) \geq 2$. In this paragraph, we explain a standard procedure that will allow us to replace $\bfG$ by its simply connected cover $\tilde\bfG$ and assume in the next paragraph that $\bfG$ is simply connected. Denote by $\tilde \bfG$ the simply connected cover of $\bfG$ and by $\tilde \pi : \tilde \bfG \to \bfG$ the corresponding central isogeny. Following \cite[Theorem I.2.3.1]{Ma91}, we have $\bfG(k)^+ = \tilde\pi(\tilde \bfG(k))$ and $\bfG(k)/\bfG(k)^+$ is a compact abelian group of finite exponent. Since $\rk(\bfG) \geq 2$, $\bfG(k)$ has property (T) and so does its lattice $\Gamma < \bfG(k)$. In particular, $\Gamma$ is finitely generated. This further implies that the image of $\Gamma$ in $\bfG(k)/\bfG(k)^+$ is finite and so $\Gamma^+ := \Gamma \cap \bfG(k)^+$ has finite index in $\Gamma$. So it is a lattice in both $\bfG(k)$ and $\bfG(k)^+$. Observe that it suffices to prove that $\bfG(k)/\bfP(k) \to  \Prob(\bfG(k)/\bfQ(k)) : x \mapsto \delta_{p(x)}$  is the essentially unique measurable $\Gamma^+$-map. Therefore, we may pull back the whole situation in $\tilde \bfG$. We then use the discussion from \cite[Example 2.14]{BBHP20}. Set $\tilde \bfP = \tilde \pi^{-1}(\bfP)$ and $\tilde \bfQ = \tilde\pi^{-1}(\bfQ)$ and $\tilde\Gamma = \tilde\pi^{-1}(\Gamma)$. Note that $\tilde\bfP < \tilde \bfQ < \tilde \bfG$ are proper parabolic $k$-subgroups and $\tilde \bfP < \tilde \bfG$ is a minimal parabolic $k$-subgroup by \cite[Theorem 22.6(i)]{Bo91}. Moreover, as $\tilde\bfG(k)$-spaces, we may identify $\tilde\bfG(k)/\tilde\bfP(k)$ with $\bfG(k)/\bfP(k) $ and $\tilde\bfG(k)/\tilde\bfQ(k) $ with $\bfG(k)/\bfQ(k)$. Since $\ker (\tilde\pi : \tilde \bfG(k) \to \bfG(k)^+)$ is finite by \cite[Corollary 2.3.2(a)]{Ma91}, it follows that $\tilde \Gamma < \tilde \bfG(k)$ is a lattice. It remains to show that the map $\tilde\bfG(k)/\tilde\bfP(k) \to \Prob(\tilde\bfG(k)/\tilde\bfQ(k)) : g\tilde\bfP(k) \mapsto \delta_{g\tilde\bfQ(k)}$ is the essentially unique measurable $\tilde\Gamma$-map. Thus, from now on and for the rest of the proof, we may assume that $\bfG = \tilde \bfG$ is simply connected.

Following \cite[Proposition 1.4]{Sh97}, we denote by $\mathbf P_0$ the $k$-{\em discompact radical} of $\mathbf P$ in the sense that $\mathbf P_0 < \mathbf P$ is the maximal $k$-subgroup with no non-trivial $k$-compact image. Moreover,  $\mathbf P_0 \lhd \mathbf P$ is a normal $k$-subgroup and $(\mathbf P/\mathbf P_0)(k)$ is compact. For notational convenience, set $G = \bfG(k)$, $Q = \bfQ(k)$, $P = \bfP(k)$ and $P_0 = \bfP_0(k)$. Denote by $p_0 : G/P_0 \to G/P : g P_0 \mapsto g P$ the canonical factor $G$-map. Since $P_0 < G$ is noncompact, Howe-Moore theorem (see \cite[Theorem 5.1]{HM77}) implies that the nonsingular action $\Gamma \curvearrowright G/P_0$ is ergodic.

Let $\Phi : G/P \to \Prob(G/Q) $ be any measurable $\Gamma$-map. Then $\Phi_0 = \Phi \circ p_0 : G/P_0 \to \Prob(G/Q)$ is a measurable $\Gamma$-map. By \cite[Theorem 4.1]{BFGW12} (whose proof carries over to the positive characteristic case thanks to \cite[Proposition 1.9 and Corollary 1.10]{BDL15}), $\Phi_0$ is a $G$-map. Since $\Phi_0 = \Phi \circ p_0$ and since $p_0 : G/P_0 \to G/P$ is a factor $G$-map, this further implies that $\Phi$ is a $G$-map. 

Fix an admissible Borel probability measure $\mu \in \Prob(G)$. Then \cite[Corollary 5.2]{BS04} implies that there exists a unique $\mu$-stationary Borel probability measure $\nu \in \Prob(G/P)$ that is $G$-quasi-invariant and such that $(G/P, \nu)$ is the $(G, \mu)$-Poisson boundary. Since any $\mu$-stationary Borel probability measure on $G/Q$ is $G$-quasi-invariant and since $G \curvearrowright G/Q$ is transitive whence ergodic, $\eta = p_\ast(\nu) \in \Prob(G/Q)$ is the unique $\mu$-stationary Borel probability measure on $G/Q$ (see e.g.\ \cite[Proposition 2.6]{BS04}). Then \cite[Corollary 2.17]{BS04} further implies that $\Phi = \{x \mapsto \delta_{p(x)}\}$. 
\end{proof}

\begin{rem}\label{rem:FT}
Let $\mathbf G$ be an almost $k$-simple connected algebraic $k$-group such that $\rk_k(\bfG) \geq 2$. Let $\Gamma < \bfG(k)$ be a lattice. We point out that in this setting, Margulis'\ factor theorem \cite[Theorem IV.2.11]{Ma91}, which is stated when $\bfG$ is simply connected, actually holds without assuming that $\bfG$ is simply connected. 

Indeed, keep the same notation as in the statement and the proof of Proposition \ref{prop-uniqueness}. Let $\rho : \bfG(k)/\bfP(k) \to X$ be a measurable $\Gamma$-map. As in the proof of Proposition \ref{prop-uniqueness}, since $\tilde \bfG(k)/\tilde \bfP(k)\cong \bfG(k)/\bfP(k)$ as $\tilde \bfG(k)$-spaces, we  can regard $\rho :  \tilde\bfG(k)/\tilde\bfP(k) \to X$ as a measurable $\tilde\Gamma$-map, where $\tilde \Gamma = \tilde \pi^{-1}(\Gamma^+)$ and $\Gamma^+ = \Gamma \cap \bfG(k)^+$. Applying Margulis' factor theorem \cite[Theorem IV.2.11]{Ma91} to the lattice $\tilde \Gamma < \tilde \bfG(k)$, there exist a parabolic $k$-subgroup $\tilde \bfP < \tilde \bfQ < \tilde \bfG$ and a $\tilde \Gamma$-equivariant measurable isomorphism $\psi : X \to \tilde \bfG(k) / \tilde \bfQ(k)$. Then we can regard $\psi : X \to \bfG(k) / \bfQ(k)$ as a $\Gamma^+$-equivariant measurable isomorphism where $\bfQ = \tilde\pi(\tilde \bfQ)$. By the proof of Proposition \ref{prop-uniqueness}, $\psi \circ \rho : \bfG(k)/\bfP(k) \to \bfG(k)/\bfQ(k)$ is the canonical factor map and so is $\Gamma$-equivariant. Since $\rho : \tilde\bfG(k)/\tilde\bfP(k) \to X$ is a factor $\Gamma$-map, it follows that $\psi : X \to \bfG(k) / \bfQ(k)$ is a $\Gamma$-equivariant measurable isomorphism.
\end{rem}

\begin{rem}\label{rem:multiplicative}
Let $\mathbf G$ be an almost $k$-simple connected algebraic $k$-group such that $\rk_k(\bfG) \geq 2$. Let $\mathbf P < \mathbf Q < \mathbf G$ be proper parabolic $k$-subgroups so that $\mathbf P < \mathbf G$ is a minimal parabolic $k$-subgroup. Let $\Gamma < \bfG(k)$ be a lattice. Set $H = \Gamma$ or $H = \bfG(k)$. Let $M$ be an ergodic $H$-von Neumann algebra and $E : M \to L^\infty(\bfG(k) / \bfP(k))$ a faithful normal ucp $H$-map. Let $\iota : L^\infty(\bfG(k) / \bfQ(k)) \hookrightarrow M$ be a $H$-equivariant unital normal embedding. Then Proposition \ref{prop-uniqueness} implies that $E \circ \iota : L^\infty(\bfG(k) / \bfQ(k)) \hookrightarrow L^\infty(\bfG(k) / \bfP(k))$ is the canonical unital normal embedding. Moreover, $L^\infty(\bfG(k) / \bfQ(k))$ lies in the multiplicative domain of $E$. Thus, upon shrinking the parabolic $k$-subgroup $\bfQ$ if necessary and using  \cite[Theorem IV.2.11]{Ma91} in case $H = \Gamma$, we may identify the multiplicative domain of $E$ with $L^\infty(\bfG(k) / \bfQ(k))$ as $H$-von Neumann algebras. Recall that the multiplicative domain of $E : M \to L^\infty(\bfG(k) / \bfP(k))$ is the largest von Neumann subalgebra of $M$ on which $E$ is multiplicative (see \cite[Definition 1.5.8]{BO08}.
\end{rem}

\subsection{The noncommutative Nevo-Zimmer theorem}

Firstly, we prove the noncommutative Nevo-Zimmer for ergodic actions of higher rank simple algebraic groups on arbitrary von Neumann algebras, namely Theorem \ref{thm:NCNZ}.

\begin{proof}[Proof of Theorem \ref{thm:NCNZ}]
We let $k$ be a local field, $\mathbf G$ be an almost $k$-simple connected algebraic $k$-group such that $\rk_k(\mathbf G) \geq 2$ and
$\mathbf P < \mathbf G$ be a minimal parabolic $k$-subgroup.
We set $G = \mathbf G(k)$ and $P = \mathbf P(k)$. 
We let $M$ be an ergodic $G$-von Neumann algebra and $E : M \to L^\infty(G/P)$ be a faithful normal ucp $G$-map.
We assume that $E$ is not $G$-invariant
and argue to show that there exist a proper parabolic $k$-subgroup $\mathbf Q$ with $\mathbf P < \mathbf Q < \mathbf G$ and a $G$-equivariant unital normal embedding $L^\infty(G/Q) \hookrightarrow M$, where $Q = \mathbf Q(k)$. By Remark \ref{rem:multiplicative}, we will obtain the last statement of the second item.

In case $\mathbf G$ is absolutely almost simple, the proof follows directly by combining 
Theorem~\ref{First half} and Theorem~\ref{thm:gauss-stat}.
Let us elaborate.
By the first theorem, we find a $G$-equivariant embedding $L^\infty(X) \to M$,
where $X$ is a $G$-space which satisfies the stabilizer assumption of the second theorem and on which the $G^+$-action is non-trivial.
Since $M$ is $G$-ergodic, we get that so is $X$.
Thus $X$ satisfies the concluded dichotomy of the second theorem.
Since the first option is ruled out by the $G^+$-non-triviality,
we get that the second option holds - 
there exist a proper parabolic $k$-subgroup ${\bfQ}<{\bfG}$
and an almost everywhere defined measurable $G$-map $X\to ({\bfG}/{\bfQ})(k)=G/Q$.
Composing the map thus obtained $L^\infty(G/Q)\to L^\infty(X)$ with the given map $L^\infty(X) \to M$,
we indeed get a $G$-equivariant unital normal embedding $L^\infty(G/Q) \hookrightarrow M$. 
This finishes the proof, under the assumption that $\mathbf G$ is absolutely almost simple.

We now turn to the general case where $\mathbf G$ is not necessarily
absolutely almost simple.
At first, we use Theorem~\ref{First half}
only to reduce the general case to the case $M=L^\infty(X, \nu)$,
where $X$ is a $G$-space on which the $G^+$-action is non-trivial.
As above, $X$ is $G$-ergodic.
We let $\bar{X}=X/\!\!/ G^+$ be the space of ergodic components of the $G^+$-action on $X$ and $\pi : X \to \bar X$ the corresponding measurable $G$-map.
We note that $\bar{X}$ is an ergodic $G/G^+$-space and $G/G^+$ is a compact group,
thus $\bar{X}$ is a transitive $G/G^+$-space.
We conclude that there exists an intermediate closed subgroup $G^+<H<G$ such that $\bar{X}\simeq G/H$ as a $G$-space. Thus, we have a measurable $G$-map $\pi : X \to G/H$. Write $X_1 = \pi^{-1}(\{H\})$ for the fiber over the base point $H \in G/H$. Observe that $X_1$ is a measurable $H$-space that is $G^+$-ergodic. Moreover, we may  regard $X = \Ind_H^G(X_1)$ as the induced $G$-space of the $H$-space $X_1$.

{\bf Claim.} There exists a faithful normal ucp $H$-map $E_1 : L^\infty(X_1) \to L^\infty(G/P)$ that is not $H$-invariant.

\begin{proof}[Proof of the Claim]
Up to taking a compact model, we may assume that $X$ is a compact metrizable $G$-space. Fix an admissible Borel probability measure $\mu \in \Prob(G)$. By \cite[Corollary 5.2]{BS04}, there exists a unique $\mu$-stationary measure $\nu_P \in \Prob(G/P)$ such that $(G/P, \nu_P)$ is the $(G, \mu)$-Poisson boundary. Set $\nu = \nu_P \circ E \in \Prob(X)$. Denote by $\beta : G/P \to \Prob(X) : b \mapsto \beta(b)$ the measurable boundary $G$-map corresponding to $E : L^\infty(X) \to L^\infty(G/P)$. By transitivity, we may assume that the map $\beta$ is everywhere defined and strictly $G$-equivariant. We have $\bary(\beta_\ast \nu_P) = \nu$. Simply denote by $m$ the unique $G$-invariant Borel probability measure on the compact $G$-space $G/H$. 

We claim that $\Map_G(G/P, \Prob(G/H)) = \{ m\}$ is the singleton that consists of the constant function equal to $m$. Indeed, let $\Phi : G/P \to \Prob(G/H)$ be any measurable $G$-map. By transitivity, we may assume that $\Phi$ is stricly $G$-equivariant and write $\eta = \Phi(P) \in \Prob(G/H)$. Then $\eta$ is $P$-invariant. Since $G^+$ acts trivially on $G/H$, $\eta$ is also $G^+$-invariant. Since $G = G^+ \cdot P$, it follows that $\eta$ is $G$-invariant and so $\eta = m$. Thus, $\Phi$ is the constant function equal to $m$. Using Furstenberg's boundary map (see e.g.\ \cite[Theorem 2.16]{BS04}), it also follows that $m \in \Prob(G/H)$ is the unique $\mu$-stationary Borel probability measure on $G/H$.

Recall that we have the measurable $G$-map $\pi : X \to G/H$. Denote by $\pi_\ast : \Prob(X) \to \Prob(G/H)$ the pushforward measurable $G$-map. Since $\pi_\ast (\nu) \in \Prob(G/H)$ is $\mu$-stationary, it follows that $\pi_\ast(\nu) = m$ by the previous paragraph. Moreover, the previous paragraph implies that the measurable $G$-map $\pi_\ast \circ \beta : G/P \to \Prob(G/H)$ is essentially constant equal to $m$. Then for almost every $b \in G/P$, we may disintegrate $\beta(b) \in \Prob(X)$ with respect to the factor map $\pi : (X, \beta(b)) \to (G/H, m)$ and write $\beta(b) = \int_{G/H} \beta(b, gH) \, {\rm d}m(gH)$ where $\beta(b, gH) \in \Prob(\pi^{-1}(gH))$. We then obtain a measurable $G$-map $\Psi : G/P \times G/H \to \Prob(X) : (b, gH) \mapsto \beta(b, gH)$. Since $H \curvearrowright G/P$ is transitive, it follows that $G \curvearrowright G/P \times G/H$ is transitive as well. We may then assume that $\Psi$ is everywhere defined and strictly $G$-equivariant. For every $gH \in G/H$, write $\nu(gH) = \int_{G/P} \beta(b, gH) \, {\rm d}\nu_P(b) \in \Prob(X)$. Observe that the everywhere defined measurable map $G/H \to \Prob(X) : gH \mapsto \nu(gH)$ gives the disintegration of $\nu \in \Prob(X)$ with respect to the factor map $\pi : (X, \nu) \to (G/H, m)$.  

Define the measurable $H$-map $\beta_1 : G/P \to \Prob(X_1) : b \mapsto \beta(b, H)$ and set $\nu_1= \bary({\beta_1}_\ast \nu_P) = \int_{G/P} \beta(b, H) \, {\rm d}\nu_P(b) \in \Prob(X_1)$. Consider the corresponding faithful normal ucp $H$-map $E_1 : L^\infty(X_1, \nu_1) \to L^\infty(G/P)$ such that $\nu_P \circ E_1 = \nu_1$. We have that $E_1$ is not $H$-invariant, because otherwise $\eta \in \Prob(X_1)$ would be $H$-invariant and thus $\nu \in \Prob(X)$ would be $G$-invariant, which in turn would imply that $E$ is $G$-invariant. We note that $E_1$ is also not $G^+$-invariant either, as $G^+ \curvearrowright G/P$ is transitive hence ergodic.
\end{proof}

We let $\tilde{\bf G}$ be the simply connected cover of ${\bf G}$ and let 
$\tilde\pi:\tilde{\bf G}\to {\bf G}$ be the corresponding central isogeny. 
We set $\tilde{G}=\tilde{\bf G}(k)$ and recall that $\tilde\pi(\tilde{G})=G^+$.
We recall that there exists a finite separable field extension $k_1$ of $k$
and an absolutely almost simple simply connected algebraic $k_1$-group ${\bf G}_1$ such that $\tilde{\bf G}$ is isomorphic to $R_k^{k_1}({\bf G}_1)$
as $k$-groups.
For a reference of this fact, see \cite[3.1.2]{Ti65}  or \cite[Ex.16.2.9]{Sp98}.
We denote by $G_1={\bf G}_1(k_1)$
and identify $\tilde{G}$ with $G_1$.
We note that the parabolic $k_1$-subgroups of ${\bf G}_1$ are in one to one correspondence with the parabolic $k$-subgroups
of $\tilde{\bf G}$ and thus also of ${\bf G}$. 
This follows for example from the discussion in \cite[\S16.2.6]{Sp98}.
In particular, if ${\bfQ}_1$ is a parabolic $k_1$-subgroup in ${\bfG}_1$, then we may identify $({\bfG}_1/{\bfQ}_1)(k_1)\simeq ({\bfG}/{\bfQ})(k)$
for some parabolic $k$-subgroup ${\bfQ}<{\bfG}$ and vice versa.
This identification is equivariant with respect to the map obtained by composing
\[ G_1 \simeq \tilde{G} \to G^+ < G. \]
In particular, we let ${\bfP}_1<{\bfG}_1$ be the minimal parabolic $k_1$-subgroup corresponding to the minimal parabolic $k$-subgroup
${\bfP}<{\bfG}$.

Next, we consider the $H$-space $X_1$ as a $G_1$-space, via the above continuous homomorphism $G_1\to H$. Set $M_1 = L^\infty(X_1)$. 
We apply Theorem~\ref{thm:gauss-stat} and we infer that there exist a proper parabolic $k_1$-subgroup $\bfP_1 < {\bfQ}_1<{\bfG}_1$ and a $G_1$-equivariant  unital normal embedding $L^\infty( ({\bfG}_1/{\bfQ}_1)(k_1)) \hookrightarrow M_1$. Then $L^\infty( ({\bfG}_1/{\bfQ}_1)(k_1))$ is necessarily contained in the multiplicative domain of 
 the faithful normal ucp $G_1$-map $E_1$ (see Remark \ref{rem:multiplicative}). We may further assume that the multiplicative domain of $E_1$ coincides with the $G_1$-von Neumann algebra $L^\infty( ({\bfG}_1/{\bfQ}_1)(k_1))$. Write ${\bfQ}<{\bfG}$ for the parabolic $k$-subgroup corresponding to the parabolic $k_1$-subgroup $\bfQ_1 < \bfG_1$. We set $Q={\bf Q}(k)$ and identify $({\bfG}/{\bfQ})(k)$ with $G/Q$ (see \cite[Proposition 20.5]{Bo91}). We may identify $L^\infty( ({\bfG}_1/{\bfQ}_1)(k_1) )$ with $L^\infty( G/Q )$ as $G_1$-von Neumann algebras. Since $E_1 : M_1 \to L^\infty( G/P)$ is $H$-equivariant, its multiplicative domain coincides with $L^\infty( G/Q )$ as $H$-von Neumann algebras. Therefore, the factor map $X_1\to G/Q$ corresponding to the embedding $L^\infty(G/Q) \hookrightarrow L^\infty(X_1)$ is $H$-equivariant.

Inducing the resulting factor $H$-map $X_1\to G/Q$ from $H$ to $G$, we get a factor $G$-map
\[ X \simeq \Ind_H^G X_1 \to \Ind_H^G G/Q \simeq G/Q \times G/H \to G/Q. \]
The proof is now complete, interpreting this factor map as a $G$-equivariant 
unital normal  embedding $L^\infty(G/Q) \hookrightarrow M$.
\end{proof}

%In the proof of the main theorem, we should pass to an abolutely simple group (we will need this later, in the algebraic part):
%Given a $k$-simple group ${\bf G}$ we can find a finite separable field extension $k_1$ of $k$
%and an absolutely almost simple $k_1$-algebraic group ${\bf G}_1$ such that ${\bf G}=R_k^{k_1}({\bf G}_1)$.
%A reference is \cite[I.1.7]{Ma91}, the end of the section.
%Probably a better reference, but I didn't check: Springer, T.A.: Linear Algebraic Groups, Ex.16.2.9.
%
%This is OK because of the following fact:
%the $k$-parabolics in $R_{k}^{k_1}({\bf G}1)$ are exactly the subgroups of the form
%$R_{k}^{k_1}({\bf Q}_1)$, where ${\bf Q}_1$ is a $k_1$-parabolic of ${\bf G}_1$.
%Here we need a reference. I found it in the following note by Brian Conrad:
%{http://virtualmath1.stanford.edu/~conrad/249BW16Page/handouts/relbruhat.pdf}
%Lemma 6.1.
%Possibly we can find it in his book "pseudo reductive groups"
%Similarly, the group $G^+$ is invariant under restriction of scalars.

%Pass to abs simple by restriction of scalars: Springer, T.A.: Linear Algebraic Groups, Ex.16.2.9. check.
%preservation of simplicity: Pi97, On Weil restriction of reductive groups and a theorem of Prasad, Theorem 1.10

Secondly, we prove the noncommutative Nevo-Zimmer for ergodic actions of lattices in higher rank simple algebraic groups on arbitrary von Neumann algebras

\begin{thm}\label{thm:NCNZ-lattices}
Let $k$ be a local field. Let $\mathbf G$ be an almost $k$-simple connected algebraic $k$-group such that $\rk_k(\mathbf G) \geq 2$ and set $G = \mathbf G(k)$. Let $\mathbf P < \mathbf G$ be a minimal parabolic $k$-subgroup and set $P = \mathbf P(k)$. Let $\Gamma < G$ be a lattice. Let $M$ be an ergodic $\Gamma$-von Neumann algebra and $E : M \to L^\infty(G/P)$ a faithful normal ucp $\Gamma$-map. The following dichotomy holds:
\begin{itemize}
\item Either $E$ is $\Gamma$-invariant.
\item Or there exist a proper parabolic $k$-subgroup $\mathbf P < \mathbf Q < \mathbf G$ and a $\Gamma$-equivariant unital normal embedding $\iota : L^\infty(G/Q) \hookrightarrow M$ where $Q = \mathbf Q(k)$ such that $E \circ \iota : L^\infty(G/Q) \hookrightarrow L^\infty(G/P)$ is the canonical unital normal embedding.
\end{itemize}
\end{thm}

\begin{proof}
We may consider the induced von Neumann algebra $G \curvearrowright \widetilde M$ together with the faithful normal ucp $G$-map $\widetilde E : \widetilde M \to L^\infty(G/P)$ and we apply Theorem \ref{thm:NCNZ}. If $\tE$ is $G$-invariant, then $E$ is $\Gamma$-invariant (see \cite[Lemma 4.6]{BBHP20}). 

Next, assume that there exist a proper parabolic $k$-subgroup $\mathbf P < \mathbf Q < \mathbf G$ and a $G$-equivariant unital normal embedding $L^\infty(G/Q) \hookrightarrow \tM$, where $Q = \mathbf Q(k)$. We note that $C(G/Q)$ is contained in the maximal compact model $\tA \subset \tM$ for the $G$-action. By Lemma \ref{MCM1}, $\tA$ is contained in $C_b(G,M)^\Gamma$. Composing the resulting embedding $C(G/Q) \to C_b(G,M)^\Gamma$ with the evaluation morphism ${\rm ev} : C_b(G,M)^\Gamma \to M$ at the identity element $1 \in G$, we obtain a unital $\ast$-homomorphism $\iota : C(G/Q) \to M$. Observe that ${\rm ev}$ is $\Gamma$-equivariant and so is $\iota$. Using Proposition \ref{prop-uniqueness}, we may extend $\iota : L^\infty(G/Q) \hookrightarrow M$ to a $\Gamma$-equivariant unital normal embedding such that $E \circ \iota : L^\infty(G/Q) \hookrightarrow L^\infty(G/P)$ is the canonical unital normal embedding. 
\end{proof}

%%%%%
\subsection{Proofs of Theorems \ref{thm:AG-simple} and \ref{thm:AG-general}}

\begin{proof}[Proof of Theorem \ref{thm:AG-simple}]
Using a combination of Theorem \ref{thm:NCNZ-lattices} and \cite[Corollary 4.20]{BBHP20}, we infer that $\Gamma$ is charmenable. Since $\Gamma$ has property (T) and since $\Rad(\Gamma) = \mathscr Z(\Gamma)$ is finite, it follows that $\Gamma$ is charfinite.
\end{proof}

\begin{proof}[Proof of Theorem \ref{thm:AG-general}]
It suffices to combine Theorem \ref{thm:AG-simple} and \cite[Theorem A]{BBHP20} following the lines of \cite[Sections 6 and 7]{BBHP20}.
\end{proof}

%%%%%
\subsection{The noncommutative Margulis factor theorem}

Connes showed in \cite{Co80} that for any icc countable discrete group $\Gamma$ with property (T), the type ${\rm II_1}$ factor $M = L(\Gamma)$ has countable outer automorphism group $\Out(M)$ and countable fundamental group $\mathscr F(M)$. This result prompted Connes to conjecture that $L(\Gamma)$ should retain $\Gamma$ (see \cite[Problem V.B.1]{Co94}).  

\begin{CRC}
Assume that $\Gamma_1$ and $\Gamma_2$ are icc countable discrete groups with property (T) such that $L(\Gamma_1) \cong L(\Gamma_2)$. Show that $\Gamma_1 \cong \Gamma_2$.
\end{CRC}

A first deep result towards Connes'\ rigidity conjecture was obtained by Cowling-Haagerup in \cite{CH88} where they showed that lattices in $\Sp(n, 1)/\{\pm 1\}$ retain the integer $n \geq 2$. In the last twenty years, there has been tremendous progress regarding the structure and the classification of group (resp.\ group measure space) von Neumann algebras thanks to Popa's deformation/rigidity theory (see the surveys \cite{Po06, Va10, Io18}). In that respect, an icc countable discrete group $\Gamma$ is said to be W$^*$-{\em superrigid} if whenever $\Lambda$ is another discrete group such that $L(\Gamma) \cong L(\Lambda)$, then $\Gamma \cong \Lambda$. Then Connes'\ rigidity conjecture asks whether every icc countable discrete group with property (T) is W$^*$-superrigid. In \cite{IPV10}, Ioana-Popa-Vaes obtained the first examples of W$^*$-superrigid groups. The examples of W$^*$-superrigid groups constructed in \cite{IPV10} are generalized wreath products and so they don't have property (T). Very recently, Chifan-Ioana-Osin-Sun constructed in \cite{CIOS21} the first class of W$^*$-superrigid groups with property (T) and thus solved Connes'\ rigidity conjecture.

Connes'\ rigidity conjecture is particularly relevant and still wide open for the class of higher rank lattices with property (T). Let $k$ be a local field. Let $\mathbf G$ be an almost $k$-simple connected algebraic $k$-group such that $\rk_k(\mathbf G) \geq 2$ and $\mathbf P < \mathbf G$ be a minimal parabolic $k$-subgroup. Set $P = \mathbf P(k)$ and $G = \mathbf G(k)$. Let $\Gamma < G$ be a lattice and denote by $\mathscr B = L( \Gamma \curvearrowright G/P)$ the group measure space von Neumann algebra associated with the nonsingular action $\Gamma \curvearrowright G/P$. In case $\mathbf G$ is center free, that is $\mathscr Z(\bfG) = \{e\}$, $L(\Gamma)$ is a type ${\rm II_1}$ factor. Moreover, the action $\Gamma \curvearrowright G/P$ is essentially free (see \cite[Lemma 6.2]{BBHP20}) and ergodic. Then $\mathscr B$ is an amenable type ${\rm III}_1$ factor (see e.g.\ the proof of \cite[Proposition 4.7]{BN11}). 

We are now ready to prove Theorem \ref{thm:NCFT}.

\begin{proof}[Proof of Theorem \ref{thm:NCFT}]
Let $L(\Gamma) \subset  M \subset \mathscr B$ be an intermediate von Neumann subalgebra. Denote by $\{u_\gamma \mid \gamma \in \Gamma\} \subset L(\Gamma)$ the canonical unitaries implementing the action $\Gamma \curvearrowright G/P$. By \cite[Lemma 3.2]{BH22}, the conjugation action $\Gamma \curvearrowright \mathscr B$ is ergodic.  Consider the canonical conditional expectation $E : \mathscr B \to L^\infty(G/P)$. Then $\Phi = E|_{ M} : M \to L^\infty(G/P)$ is a $\Gamma$-equivariant faithful normal ucp map. Using Theorem \ref{thm:NCNZ-lattices}, there are two cases to consider.

$(\rm i)$ Assume that $\Phi : M \to L^\infty(G/P)$ is invariant. Then for every $x \in M$ and every $\gamma \in \Gamma$, we have $E(x u_\gamma^*) \in \C 1$. Then \cite[Corollary 3.4]{Su18} implies that $M = L(\Gamma)$.

$(\rm ii)$ Assume that $\Phi :  M \to L^\infty(G/P)$ is not invariant. Then there exist a proper parabolic $k$-subgroup $\mathbf P < \mathbf Q < \mathbf G$ and a $\Gamma$-equivariant unital normal embedding $\iota : L^\infty(G/Q) \hookrightarrow M$ where $Q = \mathbf Q(k)$, such that $\Phi \circ \iota : L^\infty(G/Q) \hookrightarrow L^\infty(G/P)$ is the canonical unital normal embedding. Then we have $L(\Gamma \curvearrowright G/Q) \subset M \subset \mathscr B$. Since the action $\Gamma \curvearrowright G/Q$ is essentially free (see \cite[Lemma 6.2]{BBHP20}), a combination of \cite[Theorem 3.6]{Su18} and \cite[Theorem IV.2.11]{Ma91} shows that there exists a parabolic $k$-subgroup $\mathbf Q < \mathbf R < \mathbf G$ such that with $R = \mathbf R(k)$, we have $ M  = L(\Gamma \curvearrowright G/R)$. 

It is known that there are exactly $2^{\rk_k(\mathbf G)}$ intermediate parabolic $k$-subgroups $\bfP < \bfQ < \bfG$. Therefore, the $k$-rank $\rk_k(\mathbf G)$ is an invariant of the inclusion $L(\Gamma) \subset L(\Gamma \curvearrowright G/P)$.
\end{proof}

%%%%%

\bibliographystyle{plain}

\end{document}